\pdfoutput=1
\documentclass[11pt,a4paper,twoside]{article}
\usepackage{amsfonts,amsmath,amssymb}
\usepackage{graphicx}
\usepackage[font=small]{caption}
\usepackage{theorem}
\usepackage{cite}
\usepackage{indentfirst}

\usepackage[T1]{fontenc}
\usepackage[utf8]{inputenc}
\usepackage[english]{babel}
\newcommand\PaperShortTitle{Optimal recovery from inexact Fourier--Hermite coefficients}
\newif\ifsubmission
\submissiontrue
\providecommand\PaperRevision{?}
\providecommand\PaperRevisionDate{\number\year-\ifnum\month<10 0\fi\number\month-\ifnum\day<10 0\fi\number\day}
\providecommand\PaperPreparedBy{(set by make paper-pdf)}
\providecommand\PaperRef{(set by make paper-pdf)}
\InputIfFileExists{revision.tex}{}{}
\usepackage{xcolor}
\definecolor{papermetafill}{HTML}{F4F6F8}
\definecolor{papermetaframe}{HTML}{C7D0D9}
\newcommand\PaperMeta{\ifsubmission\else
  \begin{center}
  \setlength{\fboxsep}{5pt}\setlength{\fboxrule}{0.5pt}%
  \fcolorbox{papermetaframe}{papermetafill}{\footnotesize\foreignlanguage{english}{%
    \begin{tabular}{@{}l@{\hspace{10pt}}l@{}}
    \textbf{Ref} & \PaperRef\\
    \textbf{Date} & \PaperRevisionDate\\
    \textbf{Type} & working draft, edition \PaperRevision; not for distribution\\
    \textbf{Prepared} & \PaperPreparedBy
    \end{tabular}}}%
  \end{center}\fi}
\ifsubmission\else
\usepackage{fancyhdr}
\usepackage{lastpage}
\fancypagestyle{plain}{%
  \fancyhf{}%
  \fancyfoot[C]{\thepage\\[3pt]{\fontsize{6}{7}\selectfont\color{gray}%
    \PaperShortTitle\ \textbullet\ working draft, edition \PaperRevision,
    \PaperRevisionDate\ \textbullet\ \PaperPreparedBy\ \textbullet\ p.~\thepage\ of \pageref{LastPage}\\%
    \PaperRef\ \textbullet\ not for distribution}}}
\AtBeginDocument{\ifdefined\hypersetup\hypersetup{pdfauthor={\PaperPreparedBy}}%
  \else\pdfinfo{/Author (\PaperPreparedBy)}\fi}
\fi

\AtBeginDocument{%
  \setlength\abovedisplayskip{15pt plus 4pt minus 3pt}%
  \setlength\belowdisplayskip{15pt plus 4pt minus 3pt}%
  \setlength\abovedisplayshortskip{10pt plus 4pt}%
  \setlength\belowdisplayshortskip{12pt plus 4pt minus 3pt}%
}

\newenvironment{caselist}{%
  \begin{list}{}{%
    \setlength\leftmargin{1.6em}%
    \setlength\itemindent{-\leftmargin}%
    \setlength\labelwidth{0pt}%
    \setlength\labelsep{0.4em}%
    \setlength\listparindent{\parindent}%
    \setlength\topsep{\medskipamount}%
    \setlength\parsep{0pt}%
    \setlength\itemsep{\smallskipamount}%
  }}{\end{list}}

\theorembodyfont{\sl}
\newtheorem{theorem}{Theorem}[section]
\newtheorem{lemma}[theorem]{Lemma}

\newtheorem{proposition}[theorem]{Proposition}
\newtheorem{remark}[theorem]{Remark}

\numberwithin{equation}{section}

\usepackage{calc}
\numberwithin{theorem}{section}

\renewcommand{\theequation}{\arabic{section}.\arabic{equation}}

\newcommand{\card}{\operatorname{card}}
\newcommand{\bxi}{\mbox{\boldmath$\xi$\unboldmath}}

\makeatletter
\renewenvironment{thebibliography}[1]
 {\section*{\centerline{\rm\bfseries References}}%
 \@mkboth{\MakeUppercase\refname}{\MakeUppercase\refname}%
 \list{\@biblabel{\@arabic\c@enumiv}}%
 {\settowidth\labelwidth{\@biblabel{#1}}%
 \leftmargin\labelwidth
 \advance\leftmargin\labelsep
 \@openbib@code
  \usecounter{enumiv}%
  \let\p@enumiv\@empty
  \renewcommand\theenumiv{\@arabic\c@enumiv}}%
  \sloppy
  \clubpenalty4000
  \@clubpenalty \clubpenalty
  \widowpenalty4000%
  \sfcode`\.\@m
  \setlength{\itemsep}{-0.1cm}}
  {\def\@noitemerr
  {\@latex@warning{Empty 'thebibliography' environment}}%
 \endlist}
\renewcommand{\@biblabel}[1]{#1.}
\makeatother \hfuzz=0.5pt \tolerance=500

\begin{document}
\begin{center}
\textbf{Optimal recovery from inexact Fourier--Hermite\\
coefficients: minimax error and information complexity}
\end{center}
\ifsubmission\else
\begin{center}
\small Working draft, edition \PaperRevision, \PaperRevisionDate. Not for distribution.
\end{center}
\fi
\PaperMeta
\vspace*{6mm} \centerline{\textsc {S.\,G.\,Solodky $\!\!{}^{\dag}$, Y.\,A.\,Volynets $\!\!{}^{\ddag}$}}

\vspace*{5mm}
\centerline{$\!\!{}^{\dag}\!\!$ Institute of Mathematics, National Academy of Sciences of Ukraine, Kyiv}
\centerline{$\!\!{}^{\ddag}\!\!$ National University of Kyiv-Mohyla Academy, Kyiv, Ukraine}

\vspace*{8mm}
\begin{abstract}
\begin{center}
\parbox{\dimexpr\textwidth-24mm}{\small
We consider the problem of optimal recovery of functions from the weighted Wiener class $W_s^\mu$
in the uniform metric from Fourier--Hermite coefficients
$\delta$-perturbed in the $\ell_p$ norm, for $s\ne4/3$ and $p\ne4/3$.
We establish the exact order in $\delta$ of the best (minimax) recovery error.
Then, from an estimate of the minimal radius of Galerkin information, we find the order of the
information complexity, that is, of the smallest number of coefficients for which this
error is attained.
We prove that in the region $s<4/3<p$ an order-optimal method cannot be diagonal,
and that attaining the minimax order requires more coefficients
than balanced truncation uses.
We construct an explicit linear non-diagonal method that attains in this region
the minimax order of the error.
}
\end{center}
\end{abstract}

\vspace*{3mm}

\vspace*{3mm}


\vspace*{6mm}

\section{Introduction}

Suppose that the Fourier--Hermite coefficients of a function are known with an error of at most
$\delta$ in the $\ell_p$ norm. It is required to recover a continuous function so that the
largest deviation on the whole real line is as small as possible.
The properties of the function to be recovered are given by its membership in the weighted
Wiener class $W_s^\mu$. Two questions are considered. What is the smallest possible
error if all perturbed coefficients are available and arbitrary recovery operators,
including nonlinear ones, are allowed? How many coefficients are sufficient
and how many are necessary to obtain the same order of the error?

For $s\ne4/3$ and $p\ne4/3$ the answer to the first question has the form $\delta^e$
(Theorem~\ref{thm:Eopt}), and the answer to the second is $\delta^{-e/\beta}$
(Theorem~\ref{thm:complexity}); the exponents are defined in~\eqref{beta}.
These answers are linked by the estimate of the minimal radius of Galerkin information
from Theorem~\ref{thm:Rord}: when at most $N$ coefficients are used,
the smallest possible error has order $\delta^e+N^{-\beta}$.
The first term is related to the data error, the second to the restriction
on the number of coefficients used. Therefore the number of coefficients is chosen
so that the second term does not exceed the first in order.

Stable summation of orthogonal series with noisy coefficients in the uniform
metric was considered in \cite{MathePereverzev2002}, \cite{Sol&Shar}, \cite{PoPo}. The minimal radius of Galerkin
information was defined in \cite{PereverzevSolodky1996} for Fredholm equations of the first kind and
was then studied for problems of numerical differentiation and summation in
\cite{SolodkyStasyuk2024,SemenovaSolodky2026}. The smallest amount of information
that ensures a given accuracy is called, in IBC theory,
the information complexity (see, for example, \cite{Plaskota1996,SemenovaSolodky2026}).

Mean approximation by algebraic polynomials with the Chebyshev--Hermite weight and the widths
of classes of functions in the weighted spaces $L_{2,\gamma}(\mathbb R)$, $\gamma(x)=e^{-x^2}$, and
$L_{2,\gamma}(\mathbb R^2)$, $\gamma(x,y)=e^{-x^2-y^2}$, were studied in
\cite{Vakarchuk2014,Vakarchuk2022UMJ,Vakarchuk2022JMS}.
Recovery of a function from a finite set of Fourier--Hermite coefficients with random
Gaussian noise was considered in \cite{AmatoDeFeis}: there Tikhonov regularization is used,
with the parameter chosen by the GCV criterion, the function is taken from a weighted Sobolev
space with the weight $e^{-x^2}$, and what is estimated is the expectation of the squared deviation
of the regularized solution from a partial sum of the Fourier--Hermite series in the space $L_2$
with the same weight.
Exact errors and optimal
methods were found for the recovery of the identity \cite{BabenkoParfinovychSkorokhodov2021} and
diagonal \cite{BabenkoParfinovychSkorokhodov2024} operators on weighted classes of
sequences, including recovery from the first coordinates known with an error in $\ell_p$.
Uniform estimates of the Hermite functions, from which the estimates used below are derived, were obtained in
\cite{AskeyWainger1965,Muckenhoupt1970} (see also \cite[\S\,5.5]{Szego}).

The estimates of the Hermite functions, the accuracy of the methods $\mathcal S^\theta$ and the results
on diagonal methods are taken from \cite{HerMinCompanion}
(Section~\ref{imported}). In three of the four regions of the parameters $p,s$
balanced truncation attains the minimax order. The exception is
the region $s<4/3<p$: here the error can be reduced by choosing the truncation level
separately at each point $x$. An explicit linear method with such a choice is constructed
in Proposition~\ref{prop:linear3}. It recovers a continuous function
and attains the minimax order for $N\ge\delta^{-e/\beta}$.

To obtain the lower bound we use functions whose coefficients are so small
that the admissible data error can hide them completely
(Lemma~\ref{lem:rd}); this device is used in the classical theory
of optimal recovery (see, for example, \cite{MicchelliRivlin1985}).
For comparison: in the case of the Laguerre functions balanced truncation
is order-optimal for all $s$ and $p$,
if $\mu>1-1/s$ and $\theta>\mu+1/s-1$
\cite{SolodkyVolynets2026}. For the Hermite functions the contribution of the indices $k$
for which the point $x$ lies in the transition zone is essential; the corresponding function
$f_2$ is constructed in Lemma~\ref{lem:block}.

\paragraph{Terminology.}
The transition zone of the Hermite functions $h_k(x)$ is the region near the points $x^2=2k+1$.
For $x^2<2k+1$ the function $h_k(x)$ oscillates: all its zeros lie in this interval.
For $x^2>2k+1$ the oscillations stop and a fast, exponential decay of $|h_k(x)|$ begins.
In the transition zone the oscillations give way to exponential decay. For a fixed point $x$ this means
that the functions with indices $k$ noticeably smaller than $x^2/2$ are exponentially small at this point.
The estimates used are given in Lemma~\ref{lem:cap}.

The paper is organized as follows. Section~\ref{prelim} contains definitions and
auxiliary facts. Section~\ref{setting} contains the statement of the minimax problems.
Section~\ref{imported} collects auxiliary results. In Section~\ref{sharp} we obtain
estimates of the minimax error and of the minimal radius of Galerkin information, construct
an explicit linear method for region III and determine the order of the information complexity.
In Section~\ref{sec:numerics} the linear method is compared with Fourier truncation on an example from
quantum mechanics. In Appendix~\ref{app:levels} we prove Lemma~\ref{lem:levels}.

\section{Preliminaries}\label{prelim}

Denote by ${\cal L}_2(\mathbb R)$ the space of Lebesgue measurable functions $f$
on $\mathbb R$ with finite norm
\[
\|f\|_2=\left(\int_{\mathbb R}|f(x)|^2\,dx\right)^{1/2}.
\]
By ${\cal L}_\infty(\mathbb R)$ we mean the space of essentially bounded
measurable functions with the norm
\[
\|f\|_\infty=\inf\{M\ge0:\ |f(x)|\le M\ \text{almost everywhere on}\ \mathbb R\}.
\]
In both spaces functions that coincide almost everywhere are identified.
By $C(\mathbb R)$ we denote the set of real continuous functions on
$\mathbb R$. For such functions the essential supremum coincides with the supremum:
\[
\|g\|_\infty=\sup_{x\in\mathbb R}|g(x)|,
\]
and this quantity may equal $+\infty$.

The Hermite functions $h_k$ are given by the formula (see~\cite[(5.5.1), (5.5.3)]{Szego})
$$
h_k(x)=\frac{H_k(x)e^{-x^2/2}}{(2^k\,k!\sqrt\pi)^{1/2}},
\qquad
H_k(x)=(-1)^k e^{x^2}\frac{d^k}{dx^k}e^{-x^2},
\qquad k\in\mathbb N_0 ;
$$
by
$a_k(f)=\int_{\mathbb R}f(x)h_k(x)\,dx$ we denote the Fourier--Hermite coefficients of a function
$f\in {\cal L}_2(\mathbb R)$, and by $\mathbf f=(a_k(f))_{k\in\mathbb N_0}$ their
sequence. For the expansion in Hermite functions the norm $\|\cdot\|_\infty$ is the
weighted uniform metric for the Hermite polynomials with the weight $e^{-x^2/2}$.

Let $1\le s\le\infty$, $\mu>0$ and $\underline{k}=\max\{1,k\}$. We set
\[
\|f\|_{s,\mu}=\left(\sum_{k=0}^\infty
\underline{k}^{\,s\mu}|a_k(f)|^s\right)^{1/s},\qquad 1\le s<\infty,
\qquad
\|f\|_{\infty,\mu}=\sup_{k\in\mathbb N_0}\underline{k}^{\,\mu}|a_k(f)|.
\]
The corresponding space consists of the functions $f\in{\cal L}_2(\mathbb R)$ for which
this norm is finite. The class $W_s^\mu$ is its unit ball:
\[
W_s^\mu=\bigl\{f\in{\cal L}_2(\mathbb R):\ \|f\|_{s,\mu}\le1\bigr\}.
\]
Weighted Wiener spaces are defined by conditions of this type on the expansion coefficients;
their periodic version for the coefficients of the trigonometric Fourier series is considered in
\cite[\S\,2]{Kolom2023}. Here the Fourier--Hermite coefficients are used. The definition of the class coincides with
\cite[(2.6)]{HerMinCompanion}.
The conjugate exponents are $s'=s/(s-1)$ and $p'=p/(p-1)$; for $s=1$ we have
$s'=\infty$, for $s=\infty$ we have $s'=1$.

To write the estimates we need exponents related to two different
contributions of the Hermite functions. In Lemma~\ref{lem:block} these contributions give the powers
$a(s)$ and $b(s)$; the estimate over the whole class requires the larger of them,
$\bar a(s)$. The exponent $a_+(p)$ enters the estimate of the effect of the data error
under truncation of the series (statement~\textup{(I2)} below). We set
\begin{equation}\label{exps}
\begin{gathered}
a(s)=\frac34-\frac1s,
\qquad
b(s)=\frac14-\frac1{3s},
\\[4pt]
\bar a(s)=\max\{a(s),b(s)\},
\qquad
a_+(p)=\max\Bigl\{\frac34-\frac1p,\,0\Bigr\} .
\end{gathered}
\end{equation}
Since $b(s)=a(s)/3$, we have $\bar a(s)=a(s)\ge0$ for $s\ge4/3$ and $\bar a(s)=b(s)<0$ for
$s<4/3$; $a_+(p)=a(p)$ for $p>4/3$ and $a_+(p)=0$ for $p\le4/3$.

Instead of $\mathbf f$ the data $\mathbf f^\delta=\mathbf f+\bxi$ are available, where the real
error sequence $\bxi$ is not assumed to be the coefficient sequence of
any function and $\|\bxi\|_{\ell_p}\le\delta$, $1\le p\le\infty$, $0<\delta<1$.
Summation methods of the form
$$
S_N^\nu\mathbf y=\sum_{k=0}^N\nu_k^N y_k h_k
$$
form the class $\mathcal S^\theta$ if for some $C(\nu)>0$ and $\theta>0$
$$
|1-\nu_k^N|\le C(\nu)\Bigl(\frac kN\Bigr)^\theta,
\qquad 0^\theta=0,
\qquad k=0,1,\ldots,N
$$
(see, for example, \cite[(2.9), (2.10)]{HerMinCompanion}); for examples of such methods see
\cite[Example~1.1]{SolodkyVolynets2026}. A method of the class $\mathcal S^\theta$ is assumed to be
fixed in advance, and the constants depend on it only through $\theta$ and $C(\nu)$. By
$S_N$ we denote the partial sum of the Fourier--Hermite series ($\nu_k^N=1$), and $S_N f=S_N\mathbf f$.
All functions, coefficients, errors and factors are real. The symbols $c$ and $C$
denote positive constants depending only on the parameters named in the
statement containing them, and never on $N$, $\delta$, $x$, $f$, the error
sequence or the set of indices.

\subsection{Exponents of optimal recovery and the four regions}\label{sec:exps}

Under truncation of the series at the level $N$ the error estimate consists of two parts:
$N^{-\beta}$ for the discarded remainder of the series and $\delta N^{a_+(p)}$ for the error of the
retained coefficients; these estimates are given in~\textup{(I1)}, \textup{(I2)}.
Their equality determines the level $N\asymp\delta^{-1/\rho}$ and the order
$\delta^{e_{\mathrm{tr}}}$. For the minimax error we shall also need the
exponent $e$. Here and below
\begin{equation}\label{beta}
\begin{gathered}
\beta=\mu-\bar a(s),
\qquad
\rho=\beta+a_+(p),
\\[4pt]
e_{\mathrm{tr}}=\frac{\beta}{\rho},
\qquad
e=\frac{\mu-a(s)}{\mu-a(s)+a_+(p)} .
\end{gathered}
\end{equation}
In what follows we assume $\mu>\max\{0,\bar a(s)\}$.
Under this condition both quantities $\beta$ and $\rho$ are positive, and the exponents
$e_{\mathrm{tr}}$ and $e$ belong to $(0,1]$. Let us compare the exponent of balanced truncation with the exponent of the minimax
error. Their difference equals
\begin{equation}\label{eq:gap3}
e-e_{\mathrm{tr}}=\frac{\bigl(\bar a(s)-a(s)\bigr)\,a_+(p)}
{\bigl(\mu-a(s)+a_+(p)\bigr)\rho}\ge0 .
\end{equation}
The lines $s=4/3$ and $p=4/3$ divide the parameter plane into four regions: I ($s>4/3$,
$p>4/3$), II ($s>4/3$, $p<4/3$), III ($s<4/3$, $p>4/3$) and IV ($s<4/3$, $p<4/3$).
By \eqref{eq:gap3}, $e=e_{\mathrm{tr}}$ in regions I, II, IV, since there $\bar a(s)=a(s)$
or $a_+(p)=0$. In region III
$$
\bar a(s)-a(s)=b(s)-a(s)=\frac23\Bigl(\frac1s-\frac34\Bigr)>0
$$
and $a_+(p)>0$, so that $e>e_{\mathrm{tr}}$.

In what follows we shall specify a function by its coefficients. We therefore note that for
$\mu>\max\{0,\bar a(s)\}$ every sequence $a$ with $\|a\|_{s,\mu}\le1$
lies in $\ell_2$ and defines a function of the class $W_s^\mu$.
For $s\le2$ the embedding $\ell_s\subset\ell_2$ suffices.
For $2<s<\infty$ this follows from H\"older's inequality and the convergence of the series
\[
\sum_k\underline k^{-2\mu s/(s-2)}<\infty,
\qquad \mu>a(s)>\frac12-\frac1s.
\]
For $s=\infty$ we have $|a_k|\le\underline k^{-\mu}$ and
$\sum_k|a_k|^2\le\sum_k\underline k^{-2\mu}<\infty$, since $\mu>3/4$.

\section{Statement of the optimal recovery problem}\label{setting}

First we allow the whole data sequence to be used, and then
we restrict the number of its known coordinates. In both cases the result of
recovery must be a function from $C(\mathbb R)$, and the error
is measured by the norm $\|\cdot\|_\infty$.
This setting corresponds to \cite[\S\,3]{SolodkyVolynets2026}:
here the Laguerre system is replaced by the Hermite system, the space $L_{2,w}$ by
the space ${\cal L}_2(\mathbb R)$, and the norm $\|\cdot\sqrt w\|_C$ by
the norm $\|\cdot\|_\infty$. Let $1\le p\le\infty$ and $\delta>0$.

The {\it information operator} $Gf=\mathbf f$
maps ${\cal L}_2(\mathbb R)$ into $\ell_2$. Instead of the exact coefficients, the perturbed
data $\mathbf y=Gf+\bxi$ are known, where $\bxi\in\ell_p$ and $\|\bxi\|_{\ell_p}\le\delta$.
Thus the data lie in $\ell_{\max\{p,2\}}$. A {\it recovery operator} is
an arbitrary, not necessarily linear, map
$\psi:\ell_{\max\{p,2\}}\to C(\mathbb R)$.
The class of all recovery methods and their error are defined by the formulas
\[
\Psi=\bigl\{\psi G:\ \psi:\ell_{\max\{p,2\}}\to C(\mathbb R)\bigr\},
\]
\[
e_\delta(W_s^\mu,\psi G,C(\mathbb R),\ell_p)
=\sup_{f\in W_s^\mu}\ \sup_{\substack{\bxi\in\ell_p\\\|\bxi\|_{\ell_p}\le\delta}}
\bigl\|f-\psi(Gf+\bxi)\bigr\|_\infty.
\]
In the notation $\psi G$ the method is the operator $\psi$ itself on the whole data space
$\ell_{\max\{p,2\}}$, and not only the product $\psi G$ as a map on ${\cal L}_2(\mathbb R)$.
For $p\le2$ there is no difference: $G$ maps ${\cal L}_2(\mathbb R)$ onto the whole of $\ell_2$, and the product
$\psi G$ determines $\psi$ uniquely. For $p>2$ the perturbed data $Gf+\bxi$ may fail to lie
in $\ell_2$, that is, may fail to be the coefficients of any function from ${\cal L}_2(\mathbb R)$; the values
of $\psi$ on such data are not determined by the product $\psi G$, while the error $e_\delta$ depends
on them. The problem of optimal recovery from
unlimited information consists in determining the quantity
\[
\mathcal E_\delta(W_s^\mu,\Psi,C(\mathbb R),\ell_p)
=\inf_{\psi G\in\Psi}e_\delta(W_s^\mu,\psi G,C(\mathbb R),\ell_p).
\]
These definitions coincide with \cite[(3.3), (3.4)]{HerMinCompanion}.

Now let $\Omega\subset\mathbb N_0$ be a finite set of indices,
$P_\Omega$ the coordinate projection and $G_\Omega=P_\Omega G$; by $\mathbb R^\Omega$
we denote the space of tuples $(y_k)_{k\in\Omega}$ of real numbers.
From the data $G_\Omega f+P_\Omega\bxi\in\mathbb R^\Omega$ the function is recovered
by an arbitrary map $\psi:\mathbb R^\Omega\to C(\mathbb R)$.
We set
\[
\Psi_\Omega=\bigl\{\psi G_\Omega:\ \psi:\mathbb R^\Omega\to C(\mathbb R)\bigr\},
\]
\[
e_\delta(W_s^\mu,\psi G_\Omega,C(\mathbb R),\ell_p)
=\sup_{f\in W_s^\mu}\ \sup_{\substack{\bxi\in\ell_p\\\|\bxi\|_{\ell_p}\le\delta}}
\bigl\|f-\psi(G_\Omega f+P_\Omega\bxi)\bigr\|_\infty,
\qquad
\Psi_N=\bigcup_{\card(\Omega)\le N}\Psi_\Omega.
\]
Here $N\in\mathbb N$ is the maximal admissible number of coefficients used.
A method on a finite set of indices is identified with the method on the full data
sequence $\mathbf y\mapsto\psi(P_\Omega\mathbf y)$; in this sense $\Psi_N\subset\Psi$.
The problem with finite information consists in choosing both the set $\Omega$ and the recovery
operator under the constraint $\card(\Omega)\le N$.
The minimal radius of Galerkin information is
\[
R_{N,\delta}(W_s^\mu,\Psi_N,C(\mathbb R),\ell_p)=
\inf_{\psi G_\Omega\in\Psi_N}e_\delta(W_s^\mu,\psi G_\Omega,C(\mathbb R),\ell_p),
\]
and $\Psi_N\subset\Psi$ implies
\begin{equation}\label{eq:low1}
\mathcal E_\delta(W_s^\mu,\Psi,C(\mathbb R),\ell_p)\le
R_{N,\delta}(W_s^\mu,\Psi_N,C(\mathbb R),\ell_p).
\end{equation}
Now we fix the required order of the error $\delta^e$ and find
the smallest number of coefficients with which it is attainable. Let $K>0$
be fixed. The {\it information complexity} of summation is
\begin{equation}\label{eq:Nmin}
N_{\min}(\delta)=\min\left\{N\in\mathbb N:
R_{N,\delta}(W_s^\mu,\Psi_N,C(\mathbb R),\ell_p)\le K\delta^e\right\},
\end{equation}
when the set on the right is nonempty.

\section{Auxiliary results}\label{imported}

The upper bounds require estimates of individual Hermite functions
(Lemma~\ref{lem:cap}); the lower bounds require functions with coefficients on a given
finite set of indices (Lemma~\ref{lem:block}). They are followed by four
statements from \cite{HerMinCompanion}: on the remainder of the series, on the error of summation
methods, on the limits of diagonal recovery, and the block estimate of Hermite functions.
The number $m=x^2/2$ in Lemma~\ref{lem:cap} separates the indices $k$ for fixed $x$: for
$k<m$ estimate (E3) is exponentially small when $m-k$ is much larger than $m^{1/3}$, while for $k\ge m$
the power estimates (E1), (E2) hold.

\begin{lemma}\label{lem:cap}
Let $x\in\mathbb R$ and $m=x^2/2$. With absolute constants $c,C_0,C>0$ we have
\begin{itemize}
\item[\rm(E0)] $|h_k(x)|\le C_0(k+1)^{-1/12}$ for all $k\in\mathbb N_0$;
\item[\rm(E1)] $|h_k(x)|\le C(m^{1/3}+k-m)^{-1/4}$ for integers $k\ge m$;
\item[\rm(E2)] $|h_k(x)|\le C(k+1)^{-1/4}$ for integers $k\ge2m$;
\item[\rm(E3)] for $m\ge8$ and integers $0\le k<m$
$$
|h_k(x)|\le Cm^{-1/12}\exp(-cu^{3/2}),
\qquad
u=\frac{m-k}{m^{1/3}} ;
$$
in particular, $|h_k(x)|\le Ce^{-cm}$ for $k<m/2$.
\end{itemize}
\end{lemma}

\textit{Proof.} Let $\nu=2k+1$. By \cite[Lemma~4.1 (i)]{HerMinCompanion}, the
estimate $|h_k(x)|\le C(\nu^{1/3}+|x^2-\nu|)^{-1/4}$ holds for all $x$ and $k$. This implies (E0), since
$\nu^{1/3}\ge(k+1)^{1/3}$; (E1), since for $k\ge m$ we have $\nu^{1/3}\ge m^{1/3}$ and
$|x^2-\nu|=2(k-m)+1\ge k-m$; (E2), since for $k\ge2m$ we have $|x^2-\nu|\ge k+1$.

We prove (E3). Since $|h_k|$ is even, it suffices to consider $x>0$.
Set $d=m-k$; then $x^2-\nu=2d-1$. Denote $u=d/m^{1/3}$.
If $u\le1$, then $k\ge m-m^{1/3}\ge m/2$, and by (E0)
\[
|h_k(x)|\le 2^{1/12}C_0m^{-1/12}
\le 2^{1/12}\exp(1)C_0m^{-1/12}e^{-u^{3/2}}.
\]

Now let $u>1$; then $d>m^{1/3}\ge2$ and $x^2-\nu\ge d$.
We use the ready-made estimates of Askey--Wainger
\cite[table after Theorem~B, p.~700]{AskeyWainger1965}; the normalization of their Hermite functions
coincides with ours. If $r=\nu^{1/6}(x-\sqrt\nu)$, these estimates give, with absolute $A,\alpha>0$,
\[
\begin{aligned}
|h_k(x)|&\le A\nu^{-1/12}e^{-\alpha r^{3/2}},
&&\sqrt\nu\le x\le\sqrt{2\nu},\\
|h_k(x)|&\le Ae^{-\alpha x^2},
&&x\ge\sqrt{2\nu}.
\end{aligned}
\]
If \emph{$x\le\sqrt{2\nu}$,} then $m\le\nu<2m$, and
\[
r=\frac{\nu^{1/6}(x^2-\nu)}{x+\sqrt\nu}
\ge\frac{m^{1/6}d}{2\sqrt{2m}}
=\frac{u}{2\sqrt2}.
\]
Taking into account $(2\sqrt2)^{3/2}=2^{9/4}$ and $\nu^{-1/12}\le m^{-1/12}$, from the first estimate we obtain
\[
|h_k(x)|\le Am^{-1/12}e^{-c_2u^{3/2}},
\qquad
c_2=2^{-9/4}\alpha .
\]

\emph{If $x\ge\sqrt{2\nu}$,} we use the second estimate, the inequality
$u^{3/2}=d^{3/2}m^{-1/2}\le m$ and $c_2<\alpha$:
\[
|h_k(x)|\le Ae^{-2\alpha m}\le Cm^{-1/12}e^{-\alpha m}
\le Cm^{-1/12}e^{-c_2u^{3/2}}.
\]
In the case $u\le1$ the estimate above contains the factor $e^{-u^{3/2}}\le e^{-c_2u^{3/2}}$, since
$c_2<\alpha<1$. Thus (E3) is proved with the constant $c_2$.
For $k<m/2$ we have $u^{3/2}>m/2^{3/2}$; choosing $c=c_2/2^{3/2}$
and increasing $C$, we obtain both stated estimates with one constant $c$. $\Box$

\begin{figure}[htb]
\centering
\includegraphics[width=\textwidth]{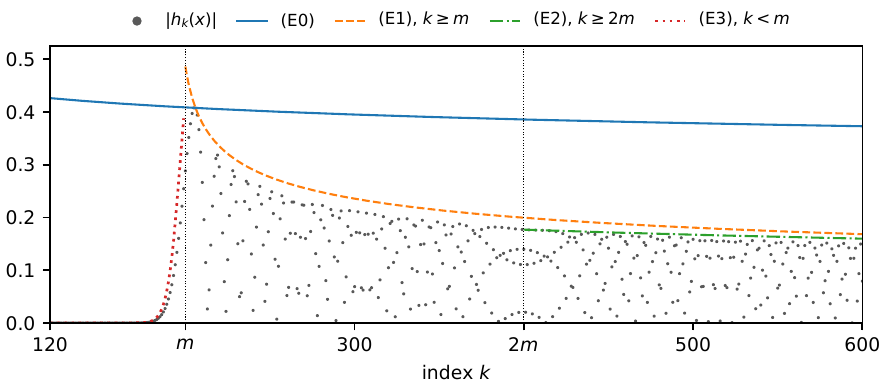}
\caption{Values of $|h_k(x)|$ for $x=20$ ($m=200$) (dots) and the right-hand sides of estimates (E0)--(E3) (lines),
each on its own range of indices $k$. The exponent $c\approx0.74$ in (E3) is fitted so that
the curve is closest to the dots for $100\le k<m$; the constants $C_0\approx0.64$ in (E0),
$C\approx0.76$ in (E1), $C\approx0.79$ in (E2) and $C\approx0.64$ in (E3) are the smallest for
which the curve is not below the dots of its range for $k\le600$. These are constants for one
value of $x$, not the absolute constants of Lemma~\ref{lem:cap}.}\label{fig:bounds}
\end{figure}

Figure~\ref{fig:bounds} shows the values of $|h_k(x)|$ at one $x$ together with the right-hand
sides of estimates (E0)--(E3).

\begin{lemma}\label{lem:block}
There exist an absolute constant $c_*>0$ and an integer $N_0$ such that for every integer
$N\ge N_0$, every set $\Omega$ of $N$ distinct integers of the interval $[N,3N]$, and all
$1\le s\le\infty$ and $\mu>0$ there exist signs $\varepsilon_k,\varepsilon'_k\in\{-1,1\}$ and
a nonempty set $W\subset\Omega$ for which the functions
$$
f_1=(3N)^{-\mu}N^{-1/s}\sum_{k\in\Omega}\varepsilon_k h_k,
\qquad
f_2=(3N)^{-\mu}(\card W)^{-1/s}\sum_{k\in W}\varepsilon'_k h_k
$$
satisfy $\|f_i\|_{s,\mu}\le1$, have Fourier--Hermite coefficients equal to zero outside $\Omega$, and
$$
\|f_1\|_\infty\ge\tfrac1{10}\,3^{-\mu}N^{-\mu+3/4-1/s},
\qquad
\|f_2\|_\infty\ge c_*\,3^{-\mu}N^{-\mu+1/4-1/(3s)} .
$$
\end{lemma}

\textit{Proof.} The quantity $N_0$, the signs $\varepsilon_k$ and the estimate for $f_1$ are taken from
\cite[Lemma~5.1 (i), (iii)]{HerMinCompanion}, where they are established for the same $N$, $\Omega$, $s$ and
$\mu>0$. Moreover, by \cite[Lemma~5.1 (ii)]{HerMinCompanion}, for the same $N$ and $\Omega$
there exist absolute constants $\kappa\in(0,1)$, $c'>0$, a set $W\subset\Omega$ with
$\card(W)\ge\kappa N^{1/3}$ and a point $x_W$ for which $|h_k(x_W)|\ge c'N^{-1/12}$ for
$k\in W$. Set $\varepsilon'_k=\operatorname{sign}h_k(x_W)$: then at the point $x_W$ all
terms of $f_2$ are positive and do not cancel each other. The coefficients of $f_2$ equal
$(3N)^{-\mu}(\card W)^{-1/s}\varepsilon'_k$ on $W$ and zero outside $W$. Since
$\underline k\le3N$ on $W$, we have $\|f_2\|_{s,\mu}\le1$, and at the point $x_W$
$$
\|f_2\|_\infty\ge(3N)^{-\mu}(\card W)^{-1/s}\sum_{k\in W}|h_k(x_W)|
\ge c'\,3^{-\mu}N^{-\mu-1/12}(\card W)^{1-1/s} .
$$
Since $1-1/s\ge0$, $\kappa<1$ and $\card(W)\ge\kappa N^{1/3}$, this is at least
$c'\kappa\,3^{-\mu}N^{-\mu+1/4-1/(3s)}$. $\Box$

Note that the function $f_1$ will give the lower bound in $\delta$ in region III (Theorem~\ref{thm:Eopt}), and
the lower bound in $N$ (Theorem~\ref{thm:Rord}) will be given by $f_1$ for $s>4/3$ and by $f_2$ for $s<4/3$.

Four more statements from \cite{HerMinCompanion} are used; in them $s\ne4/3$, $p\ne4/3$,
$\mu>\max\{0,\bar a(s)\}$ and $0<\delta<1$.

\begin{itemize}
\item[(I1)] \cite[Lemma~4.4]{HerMinCompanion}: for $f\in W_s^\mu$ the Fourier--Hermite series of $f$
converges absolutely and uniformly on $\mathbb R$, its sum is continuous, and
\[\|f-S_N f\|_\infty\le c\|f\|_{s,\mu}N^{-\beta}\]
for $N\ge1$, where $c=c(s,\mu)$, $\beta=\mu-\bar a(s)$ (see~\eqref{beta}).

\item[(I2)] \cite[Theorem~4.7]{HerMinCompanion}: every method of the class $\mathcal S^\theta$
with $\theta>\beta$ satisfies
\[\|f-S_N^\nu\mathbf f^\delta\|_\infty\le c\,(N^{-\beta}+\delta N^{a_+(p)})\]
for
all $f\in W_s^\mu$, $\|\bxi\|_{\ell_p}\le\delta$ and $N\ge1$, and for
$\gamma_1\delta^{-1/\rho}\le N\le\gamma_2\delta^{-1/\rho}$, $0<\gamma_1\le\gamma_2$, it satisfies the estimate
\[\|f-S_N^\nu\mathbf f^\delta\|_\infty\le c\,\delta^{e_{\mathrm{tr}}}, \quad
e_{\mathrm{tr}}=\beta/\rho;\]
the constants depend on
$s$, $p$, $\mu$, $\theta$, $C(\nu)$, and in the second estimate also on $\gamma_1$, $\gamma_2$.
We call such a method with level $N=\lceil\delta^{-1/\rho}\rceil$ \emph{balanced truncation};
it uses $\lceil\delta^{-1/\rho}\rceil+1\le3\delta^{-1/\rho}$ Fourier--Hermite coefficients.

\item[(I3)] \cite[Theorem~3.1]{HerMinCompanion}: the error of every diagonal method
\[\psi_\nu(\mathbf y)=\sum_{k\in\Omega}\nu_k y_k h_k\]
with a finite nonempty
$\Omega\subset\mathbb N_0$ and real $\nu_k$ is at least $c\,\delta^{e_{\mathrm{tr}}}$, where
$c=c(s,p,\mu)$ does not depend on $\Omega$ and $\nu$. For $s>4/3$ or $p<4/3$ we have
\[\mathcal E_\delta(W_s^\mu,\Psi,C(\mathbb R),\ell_p)\asymp\delta^{e_{\mathrm{tr}}},\]
and this order is attained by balanced truncation.

\item[(I4)] \cite[Lemma~4.3]{HerMinCompanion} and the scheme of the proof of
\cite[Lemma~4.4]{HerMinCompanion}: for all $x\in\mathbb R$ and integers $N\ge1$
\[
\Bigl(\sum_{k\ge N}k^{-\mu s'}|h_k(x)|^{s'}\Bigr)^{1/s'}\le CN^{-\beta},
\qquad C=C(s,\mu),
\]
where for $s=1$ the left-hand side means $\sup_{k\ge N}k^{-\mu}|h_k(x)|$. Moreover,
by \cite[Lemma~4.3]{HerMinCompanion}, for all integers $M\ge1$
\[
\Bigl(\sum_{M\le k<2M}|h_k(x)|^{s'}\Bigr)^{1/s'}\le CM^{\bar a(s)}
\]
(for $s=1$, the maximum). On the block $E_i=\{k:\ 2^{i}N\le k<2^{i+1}N\}$ the factor
$k^{-\mu}$ is at most $(2^{i}N)^{-\mu}$. Therefore the terms of the left-hand side of (I4) with indices
$k\in E_i$ give
\[
\Bigl(\sum_{k\in E_i}k^{-\mu s'}|h_k(x)|^{s'}\Bigr)^{1/s'}
\le C(2^{i}N)^{-\mu+\bar a(s)}=C(2^{i}N)^{-\beta}.
\]
The blocks $E_i$, $i\ge0$, cover all indices $k\ge N$, and by the triangle inequality in
$\ell_{s'}$ the left-hand side of (I4) is at most the sum of these estimates over $i\ge0$. Since $\beta>0$,
this sum is a convergent geometric progression with ratio $2^{-\beta}<1$:
\[
\sum_{i\ge0}C(2^{i}N)^{-\beta}=C\bigl(1-2^{-\beta}\bigr)^{-1}N^{-\beta}.
\]
\end{itemize}

\section{Lower and order-optimal bounds for the minimax quantities}\label{sharp}

Throughout this section $s\ne4/3$, $p\ne4/3$, $\mu>\max\{0,\bar a(s)\}$, the number $N\ge1$ is an integer,
$0<\delta<1$ and $N_0$ is from Lemma~\ref{lem:block}. In region III a threshold
level $\delta_0\in(0,1)$ is introduced below: the estimates are proved first for $\delta<\delta_0$, then for
$\delta_0\le\delta<1$.
The constants $c$, $C$ depend only on $s$, $p$, $\mu$ unless stated otherwise.

\subsection{Unlimited information}\label{sec:unlim}

First we examine which functions cannot be distinguished by the data.
If $\|G_\Omega g\|_{\ell_p}\le\delta$, then the admissible error
$-G_\Omega g$ turns the observed coefficients of $g$ into the sequence
of all zeros. The same data are obtained for $-g$ with the opposite
error. Therefore the norm $\|g\|_\infty$ of such functions bounds from below the error of
every method. For a finite set $\Omega\subset\mathbb N_0$ or
$\Omega=\mathbb N_0$ (then $G_\Omega=G$) we set

\begin{equation}\label{diam}
\begin{gathered}
D_\Omega(\delta)=\sup\bigl\{\|g\|_\infty:\ g\in W_s^\mu,\ \|G_\Omega g\|_{\ell_p}\le\delta\bigr\},
\\[4pt]
D(\delta)=D_{\mathbb N_0}(\delta),
\qquad
B=\sup_{f\in W_s^\mu}\|f\|_\infty ;
\end{gathered}
\end{equation}
here $\|\cdot\|_{\ell_p}$ of a sequence may equal $+\infty$. The quantity $B$ does not
depend on $\delta$ and $N$; its estimate will be needed for large $\delta$ and for
$N=1$. The set
\begin{equation}\label{consist}
A_\Omega(y)=\bigl\{f\in W_s^\mu:\ \|G_\Omega f-y\|_{\ell_p}\le\delta\bigr\},
\qquad A(y)=A_{\mathbb N_0}(y),
\end{equation}
consists of the functions consistent with the data $y$. An operator that selects an element of
$A_\Omega(y)$ when this set is nonempty and the zero function otherwise is called a {\it consistent
recovery operator}.
It is needed for the upper bounds: by Lemma~\ref{lem:rd}(b) its error is at most
$2D_\Omega(\delta)$, and by Lemma~\ref{lem:rd}(a) the error of every method is at least
$D_\Omega(\delta)$, so that from now on it suffices to estimate $D_\Omega(\delta)$. Part~(c) of the same
lemma compares $D_{\Omega_N}$ for the first $N$ indices with $D$ and is used in estimating the
radius of Galerkin information.

\begin{lemma}\label{lem:rd}
\begin{itemize}
\item[\rm(a)] For every recovery operator $\psi$
$$
\sup\bigl\{\|f-\psi(G_\Omega f+P_\Omega\bxi)\|_\infty:\ f\in W_s^\mu,\
\|\bxi\|_{\ell_p}\le\delta\bigr\}\ge D_\Omega(\delta)
$$
{\rm(cf.~\cite[Theorem~2.1]{MicchelliRivlin1985})}.
\item[\rm(b)] If $f\in W_s^\mu$, $\|\bxi\|_{\ell_p}\le\delta$ and
$u\in A_\Omega(G_\Omega f+P_\Omega\bxi)$, then $\|f-u\|_\infty\le2D_\Omega(\delta)$.
\item[\rm(c)] $B<\infty$, and for $\Omega_N=\{0,1,\ldots,N-1\}$, $N\ge2$,
$$
D_{\Omega_N}(\delta)\le D(\delta)+cN^{-\beta},
$$
where $c$ depends only on $s$ and $\mu$.
\end{itemize}
\end{lemma}

\textit{Proof.} The class $W_s^\mu$ is convex and symmetric, and the functional
$\|\cdot\|_{\ell_p}$ with values in $[0,\infty]$ satisfies the triangle inequality.

(a) Let $g\in W_s^\mu$ and $\|G_\Omega g\|_{\ell_p}\le\delta$. The errors
$\bxi=\mp G_\Omega g$ are admissible, the functions $\pm g$ give the same data $\mathbf 0$, and for
$u=\psi(\mathbf 0)$ we have
\[\max\{\|g-u\|_\infty,\|g+u\|_\infty\}\ge\|g\|_\infty.\]

(b) Let $y=G_\Omega f+P_\Omega\bxi$. The condition $u\in A_\Omega(y)$ and definition
\eqref{consist} imply that $u\in W_s^\mu$. Therefore
$-u\in W_s^\mu$, and hence the function $g=(f-u)/2$ also lies in $W_s^\mu$. Consequently,
\[\|G_\Omega g\|_{\ell_p}\le\frac12(\|P_\Omega\bxi\|_{\ell_p}+\|G_\Omega u-y\|_{\ell_p})\le\delta,\]
whence $\|f-u\|_\infty=2\|g\|_\infty\le2D_\Omega(\delta)$.

(c) For $f\in W_s^\mu$ we have $|a_0(f)|,|a_1(f)|\le1$, therefore by (E0) and (I1) with $N=1$ we get
$\|f\|_\infty\le2C_0+c_1$, where $c_1$ is the constant of (I1); hence $B\le2C_0+c_1<\infty$.

Let $g\in W_s^\mu$ and $\|G_{\Omega_N}g\|_{\ell_p}\le\delta$. The coefficients of the function
$g_N=S_{N-1}g$ equal $a_k(g)$ for $k\le N-1$ and zero for $k\ge N$. Therefore
$\|g_N\|_{s,\mu}\le\|g\|_{s,\mu}\le1$ and $\|Gg_N\|_{\ell_p}=\|G_{\Omega_N}g\|_{\ell_p}\le\delta$,
so that $g_N$ is admissible in the definition of $D(\delta)$ and $\|g_N\|_\infty\le D(\delta)$.
Since $N\ge2$, we have $N-1\ge N/2$, and by (I1)
\[\|g-g_N\|_\infty\le c_1(N-1)^{-\beta}\le2^\beta c_1N^{-\beta}.\]
Therefore
\[\|g\|_\infty\le\|g_N\|_\infty+\|g-g_N\|_\infty\le D(\delta)+2^\beta c_1N^{-\beta}.\]
The supremum over such $g$ gives (c) with $c=2^\beta c_1$. $\Box$

\smallskip
\noindent\textit{Region III.} Let $s<4/3<p$, so that $s'>4$, $1\le p'<4$ and
$a(s)<0<a(p)$.
We estimate $D(\delta)$. For an admissible function $g$ both the weighted
$\ell_s$ norm of the coefficients and their $\ell_p$ norm are bounded. We apply the first bound
to the remainder of its series, and the second to its beginning. The splitting level $L\ge-1$ is a real
number, not necessarily an integer; we will choose it separately at each point $x$.
For these two parts we introduce the quantities
\[
\mathcal T(x,L)=\Bigl(\sum_{k>L}\underline k^{\,-\mu s'}|h_k(x)|^{s'}\Bigr)^{1/s'},
\qquad
\mathcal H(x,L)=\Bigl(\sum_{0\le k\le L}|h_k(x)|^{p'}\Bigr)^{1/p'},
\qquad L\ge-1.
\]
For $L=-1$ we have $\mathcal H(x,-1)=0$, and $\mathcal T(x,-1)$ involves
all indices $k\ge0$. For $s=1$ the quantity $\mathcal T$ means
$\sup_{k>L}\underline k^{-\mu}|h_k(x)|$, and for $p=\infty$
the quantity $\mathcal H$ equals $\sum_{k\le L}|h_k(x)|$.

Let $\bigl\|(\underline k^{\,\mu}a_k)\bigr\|_{\ell_s}\le1$ and
$\|\bxi\|_{\ell_p}\le\delta$. H\"older's inequality with exponents $s$, $s'$ in the remainder of the series
and $p$, $p'$ in its beginning gives, for every real $L\ge-1$,
\begin{equation}\label{split}
\Bigl|\sum_{k>L}a_k h_k(x)+\sum_{0\le k\le L}\xi_k h_k(x)\Bigr|
\le\mathcal T(x,L)+\delta\,\mathcal H(x,L).
\end{equation}
The Fourier--Hermite series of an admissible function $g$ converges absolutely, and its continuous sum
represents $g$ (I1); we denote its value at the point $x$ by $g(x)$. Therefore \eqref{split} with $a_k=\xi_k=a_k(g)$ gives
\begin{equation}\label{holder}
|g(x)|\le\mathcal T(x,L)+\delta\,\mathcal H(x,L).
\end{equation}
We call the first term the estimate of the remainder of the series, and the second the error estimate: it
accounts for the contribution of the data error $\bxi$, $\|\bxi\|_{\ell_p}\le\delta$, to the
coefficients with $k\le L$. As $L$ increases, the first estimate does not increase, and the second
does not decrease. For sufficiently small $\delta$ we choose the level $L$ so that the right-hand side is
at most $C\delta^e$ for all $x$; then the same estimate holds for $D(\delta)$.
Indeed, the level $L$ depends only on $x$ and $\delta$, and the constant $C$ depends
neither on $x$ nor on $g$. Therefore $|g(x)|\le C\delta^e$ for every admissible function $g$ and
every point $x$; the supremum over $x$ gives $\|g\|_\infty\le C\delta^e$, and the supremum over
admissible $g$ in~\eqref{diam} gives $D(\delta)\le C\delta^e$.
Here $L$ is an auxiliary truncation boundary; the number of data used will be
bounded separately.

\emph{How the level will be chosen.} We now give a preliminary explanation of
how $L$ is chosen; the choice itself will be justified in Lemma~\ref{lem:levels}.
Set $m=x^2/2$. By (E3) the values $|h_k(x)|$ with $k<m$ are exponentially small when
$m-k$ is much larger than $m^{1/3}$, and by (E1) for $k\ge m$ they are at most
$C(m^{1/3}+k-m)^{-1/4}$. Therefore it is convenient to write the level as $L=m+t$, $t>0$. Let $m\ge8$
and $m^{1/3}\le t\le m$. In the remainder of the series $\underline k^{\,-\mu}\le m^{-\mu}$, and the sum
$\sum_{k>m+t}(k-m)^{-s'/4}$ for $s'>4$ is of order $t^{1-s'/4}$; in the beginning of the series
the sum $\sum_{m\le k\le m+t}(m^{1/3}+k-m)^{-p'/4}$ for $p'<4$ is of order
$t^{1-p'/4}$, while the indices $k<m$ give at most the same order. Since $1/s'-1/4=a(s)$ and
$1/p'-1/4=a(p)$, one should expect the estimates
\begin{equation}   \label{2est}
\mathcal T(x,m+t)\le Cm^{-\mu}t^{a(s)},
\qquad
\delta\,\mathcal H(x,m+t)\le C\delta\,t^{a(p)} .
\end{equation}
This is a heuristic; these estimates will be obtained rigorously in the proof of Lemma~\ref{lem:levels}.

We equate the right-hand sides of estimates (\ref{2est}) without the constants.
The balance point $t_*(m)$ is determined by the equality
$m^{-\mu}t^{a(s)}=\delta t^{a(p)}$. We denote the difference of the exponents of the two estimates by
\[
\lambda=a(p)-a(s)=\frac1s-\frac1p>0 ;
\]
then this equality is equivalent to $t^{\lambda}=m^{-\mu}/\delta$, and $t_*(m)$ decreases in $m$.
Denote by $M_0$ and $M_1$
the values of $m$ at which $t_*(m)$ reaches the bounds
$m$ and $m^{1/3}$, respectively:
\begin{equation}\label{M0}
\begin{gathered}
t_*(m)=\Bigl(\frac{m^{-\mu}}{\delta}\Bigr)^{1/\lambda},\qquad m>0,
\\[4pt]
M_0=\delta^{-1/(\mu+\lambda)},
\qquad
M_1=\delta^{-1/(\mu+\lambda/3)}.
\end{gathered}
\end{equation}
Thus $t_*(M_0)=M_0$, $t_*(M_1)=M_1^{1/3}$, and $m^{1/3}\le t_*(m)\le m$
exactly for $M_0\le m\le M_1$.
At the first of these levels both estimates in (\ref{2est}) equal $\delta^e$. In region III $a_+(p)=a(p)$,
so the denominator of the exponent $e$ in~\eqref{beta} equals $\mu-a(s)+a(p)=\mu+\lambda$, and
\begin{equation}\label{M0e}
e=\frac{\mu-a(s)}{\mu+\lambda},
\qquad
M_0^{a(s)-\mu}=\delta M_0^{a(p)}=\delta^{e} .
\end{equation}
Here the first equality for $M_0$ is the balance equality at $m=t=M_0$ (since
$t_*(M_0)=M_0$), and the second is obtained by substituting $M_0=\delta^{-1/(\mu+\lambda)}$ into
$M_0^{a(s)-\mu}$.
We choose the truncation boundary over three ranges of $m$:
\[
\Lambda(m)=
\begin{cases}
2M_0, & 0\le m\le M_0,\\
m+t_*(m), & M_0<m\le M_1,\\
m+m^{1/3}, & m>M_1.
\end{cases}
\]
In the middle range the balance point is used. For $m\le M_0$ we have $t_*(m)\ge m$,
and the constant level $2M_0=M_0+t_*(M_0)$ is taken; since $2M_0\ge2m$, the remainder of the series
contains only indices $k>2m$, where (E2) holds. For $m>M_1$ we have $t_*(m)<m^{1/3}$, and
the level is kept at the boundary of the transition zone.

In Lemma~\ref{lem:levels} the estimate at the level $\Lambda(m)$ is proved for $m<2M_1$, although for
$m>M_1$ the whole series is already small: the overlap $M_1<m<2M_1$ will allow the explicit method of
Proposition~\ref{prop:linear3} to pass continuously from the partial sum to zero.

\emph{Threshold level $\delta_0$.} Set
\[
\delta_0=\max\{N_0,8\}^{-(\mu+\lambda)}\in(0,1).
\]
Since $M_0=\delta^{-1/(\mu+\lambda)}$, for $\delta<\delta_0$ we have $M_0>\max\{N_0,8\}$;
this is needed to apply Lemma~\ref{lem:block} and the estimates of Lemma~\ref{lem:cap} for $m\ge8$.

\emph{Large $\delta$.} For $\delta_0\le\delta<1$ the order $\delta^e$ follows from rough
estimates. If the quantity being bounded from above (for example, $D(\delta)$ or the error of a method)
is, for such $\delta$, at most a constant $C$, then, since $e>0$ and
$\delta/\delta_0\ge1$,
\[
C\le C\delta_0^{-e}\delta^{e}.
\]
It suffices to prove the lower bound $\mathcal E_\delta\ge c\,\delta^e$ for $\delta<\delta_0$.
Indeed, the error of every method does not decrease in $\delta$, since the supremum is taken
over a larger ball of errors; therefore $\mathcal E_\delta$ does not decrease either, and for
$\delta_0\le\delta<1$
\[
\mathcal E_\delta\ge\mathcal E_{\delta_0/2}\ge c\,(\delta_0/2)^{e}\ge c\,(\delta_0/2)^{e}\delta^{e}.
\]
The constants here depend on $\delta_0$, that is, only on $s$, $p$, $\mu$, since $N_0$ is absolute.

\begin{lemma}\label{lem:levels}
Let the parameters $p, s$ lie in region III, $x\in\mathbb R$, $m=x^2/2$ and $\delta<\delta_0$.
Then $8<M_0<M_1$, and
\begin{equation}\label{eq:lin-P1}
\mathcal T(x,\Lambda(m))+\delta\,\mathcal H(x,\Lambda(m))\le C\delta^{e},
\qquad m<2M_1,
\end{equation}
\begin{equation}\label{eq:lin-P2}
\mathcal T(x,-1)\le C\delta^{e},
\qquad m>M_1 .
\end{equation}
\end{lemma}

The proof of the lemma is given in Appendix~\ref{app:levels}.

\vskip 2mm

\begin{theorem}\label{thm:Eopt}
In all four regions
\begin{equation}\label{reduction}
D(\delta)\le\mathcal E_\delta(W_s^\mu,\Psi,C(\mathbb R),\ell_p)\le2D(\delta),
\qquad
c\,\delta^{e}\le\mathcal E_\delta(W_s^\mu,\Psi,C(\mathbb R),\ell_p)\le C\,\delta^{e} .
\end{equation}
The upper bounds are provided by any consistent recovery operator with $\Omega=\mathbb N_0$.
In particular, such an operator is given by the minimization problem
\begin{equation}\label{convex}
\widehat a(y)=\arg\min\bigl\{\|a\|_{s,\mu}:\ \|a-y\|_{\ell_p}\le\delta\bigr\}.
\end{equation}
The function is recovered from the sequence $\widehat a(y)$ thus found (the minimum is attained when $A(y)$
is nonempty; otherwise the zero function is recovered). For $s>4/3$ or $p<4/3$ the same order is given by the method (I2) at
every level $\gamma_1\delta^{-1/\rho}\le N\le\gamma_2\delta^{-1/\rho}$ with
$0<\gamma_1\le\gamma_2$, in particular by balanced truncation.
\end{theorem}

\textit{Proof.} \emph{Bounds via the quantity $D(\delta)$.} The lower bound
in the left inequality \eqref{reduction} follows from Lemma~\ref{lem:rd}(a) with $\Omega=\mathbb N_0$.
Next, let $\psi$ be a consistent recovery operator with
$\Omega=\mathbb N_0$; its values lie in $W_s^\mu$ by the definition of a consistent operator and are continuous by (I1), so that $\psi G\in\Psi$. On admissible
data $y=Gf+\bxi$, Lemma~\ref{lem:rd}(b) gives the upper bound in the left relation \eqref{reduction}.

\emph{Lower bound in $\delta$.} In regions I, II and IV it follows from (I3), since there
$e=e_{\mathrm{tr}}$ (Section~\ref{sec:exps}). In region III we take the function $f_1$
of Lemma~\ref{lem:block} with indices of order $M_0$: the $\ell_p$ norm of its coefficients is at most
$\delta$, so the data error can hide them completely, while $\|f_1\|_\infty\ge c\,\delta^e$. Let first $\delta<\delta_0$ and
$N=\lceil M_0\rceil$. Then $M_0>N_0$, so $N\ge N_0$ and $M_0\le N\le2M_0$. We apply
Lemma~\ref{lem:block} to $\Omega=\{N,\ldots,2N-1\}$. The function $f_1$ lies in $W_s^\mu$, and its
coefficients equal $\pm(3N)^{-\mu}N^{-1/s}$ on $N$ indices, therefore
$$
\|\mathbf f_1\|_{\ell_p}=3^{-\mu}N^{-(\mu+\lambda)}\le M_0^{-(\mu+\lambda)}=\delta .
$$
Hence $f_1$ is admissible in the definition of $D(\delta)$. Since
$a(s)-\mu<0$, by \eqref{M0e}
$$
\mathcal E_\delta\ge D(\delta)\ge\|f_1\|_\infty\ge\tfrac1{10}\,3^{-\mu}N^{a(s)-\mu}
\ge\tfrac1{10}\,3^{-\mu}2^{a(s)-\mu}M_0^{a(s)-\mu}=c\,\delta^{e} .
$$

\emph{Upper bound in $\delta$ outside region III.}
In regions I, II and IV $e=e_{\mathrm{tr}}$, and by (I3) and
(I2) $\mathcal E_\delta\le C\delta^{e}$, while the method (I2) at the indicated levels gives the same
order. Let $\psi$ be a consistent recovery operator with $\Omega=\mathbb N_0$.
By Lemma~\ref{lem:rd}(b) its error is at most $2D(\delta)$, and by the left inequality
\eqref{reduction} $D(\delta)\le\mathcal E_\delta$; therefore
\[e_\delta(W_s^\mu,\psi G,C(\mathbb R),\ell_p)\le2D(\delta)\le2\mathcal E_\delta\le C\delta^{e}.\]

\emph{Upper bound in $\delta$ in region III.}
Let first $\delta<\delta_0$. By \eqref{reduction} it suffices
to prove that $D(\delta)\le C\delta^{e}$.
Let $g$ be admissible in the definition of $D(\delta)$, $x\in\mathbb R$ and $m=x^2/2$. Estimate
\eqref{holder} holds for $g$ at every level $L\ge-1$. We choose the level $L$ according to $m$:
for $m\le M_1$ we take $L=\Lambda(m)$ and apply \eqref{eq:lin-P1}, for $m>M_1$ we take $L=-1$ and
apply \eqref{eq:lin-P2}. Hence $|g(x)|\le C\delta^{e}$ for all $x$.

Now let $\delta_0\le\delta<1$ in region III\@. We apply the paragraph ``Large $\delta$'' of
Section~\ref{sec:unlim}. The lower bound $\mathcal E_\delta\ge c\,\delta^e$ is proved above for
$\delta_0/2<\delta_0$, therefore it also holds for $\delta\ge\delta_0$ with a constant
depending on $\delta_0$. For the upper bound a constant suffices ($B<\infty$ by
Lemma~\ref{lem:rd}(c)). Let $\psi$ be a consistent recovery operator with
$\Omega=\mathbb N_0$. By Lemma~\ref{lem:rd}(b)
\[
e_\delta(W_s^\mu,\psi G,C(\mathbb R),\ell_p)\le2D(\delta)\le2B,
\qquad
\mathcal E_\delta\le e_\delta(W_s^\mu,\psi G,C(\mathbb R),\ell_p).
\]
Therefore \eqref{reduction} and the attainment of the upper bound by every consistent recovery
operator hold for all $0<\delta<1$, and the constants depend only on $s$, $p$,
$\mu$, since $N_0$ is absolute.

\emph{Choice of the function from the data.} The upper bounds are proved above for every
consistent recovery operator. It remains to verify the claim of the theorem that
such an operator is given by problem \eqref{convex}. For this we must show that for
nonempty $A(y)$ the minimum in \eqref{convex} is attained and that the sequence found,
$\widehat a(y)$, consists of the coefficients of a function lying in $A(y)$. The main part is
the attainment of the minimum: there are infinitely many coordinates, so this requires proof.

The idea is as follows. Take sequences admissible in \eqref{convex} whose norms
approach the infimum. Their coordinates are bounded, so one can choose a
subsequence in which every coordinate converges. Both constraints of the problem ---
on the norm $\|a\|_{s,\mu}$ and on the distance $\|a-y\|_{\ell_p}$ --- are sums or suprema of
nonnegative quantities, each of which depends on one coordinate. Such
constraints are preserved under the coordinatewise passage to the limit. Therefore the limit is admissible,
and its norm is at most the infimum.

In more detail. Let $v$ be the infimum in \eqref{convex}. The coefficients of any $f\in A(y)$
are admissible, therefore $v\le1$. Take admissible $a^{(j)}$ with
$\|a^{(j)}\|_{s,\mu}\to v$ and $\|a^{(j)}\|_{s,\mu}<2$; then
$|a^{(j)}_k|\le\underline k^{\,\mu}|a^{(j)}_k|<2$ for all $k$ and $j$. By the
Bolzano--Weierstrass theorem and the diagonal process there is a subsequence in which
every coordinate converges; let $a_k$ be the limits. For $s<\infty$ and every integer $J$
the sum $\sum_{k\le J}\underline k^{\,s\mu}|a_k|^s$ is the limit of the same sums for $a^{(j)}$
and is therefore at most $v^s$; letting $J\to\infty$ we obtain $\|a\|_{s,\mu}\le v$. Likewise
$\|a-y\|_{\ell_p}\le\delta$. For $s=\infty$ or $p=\infty$ we pass to the limit in each
coordinate and take the supremum over $k$. Hence $a$ is admissible and $\|a\|_{s,\mu}=v$: the minimum
is attained. Since $v\le1$, by the remark at the end of Section~\ref{sec:exps} the sequence $a$
defines a function $u\in W_s^\mu$; its series converges to a continuous sum by (I1), and $u\in A(y)$.
For empty $A(y)$ we set $u=0$. If there are several minimizers, any of them is taken.
The map $y\mapsto u$ is a consistent recovery operator, and the upper bound has already been proved for it. $\Box$

\vskip 2mm

\subsection{Radius of Galerkin information}

When we pass to a finite set of coefficients, one more restriction is added
to the data error: the coefficients with unused indices
remain unknown even in the absence of noise. We show that the corresponding
error is of order $N^{-\beta}$ and adds up in order
with the already found quantity $\delta^e$.

\begin{theorem}\label{thm:Rord}
In all four regions of the values of the parameters $p, s$,
\begin{equation}\label{eq:Rord}
c\left(\delta^{e}+N^{-\beta}\right)\le R_{N,\delta}(W_s^\mu,\Psi_N,C(\mathbb R),\ell_p)
\le C\left(\delta^{e}+N^{-\beta}\right),
\end{equation}
and the upper bound is delivered by any consistent recovery operator with $\Omega=\Omega_N=\{0,1,\ldots,N-1\}$;
it lies in $\Psi_N$.
\end{theorem}

\textit{Proof.} The bound $R_{N,\delta}\ge c\,\delta^{e}$ follows from \eqref{eq:low1} and
Theorem~\ref{thm:Eopt}. For the bound $R_{N,\delta}\ge c\,N^{-\beta}$ we construct
a function all of whose observed coefficients vanish, while the function itself is sufficiently large
(cf.~\cite[Theorem~4.10]{SolodkyVolynets2026}, where the lower bound for Laguerre functions is constructed in this way).
Even exact observations do not distinguish this function from its opposite.
Assume first that
$N\ge N_0$, $\card(\Omega)\le N$ and $\psi G_{\Omega}\in\Psi_{\Omega}$. The interval
$[N,3N]$ contains $2N+1$ integer points, and after removing the points of $\Omega$ one can
choose $N$ of them; denote their set by $\Omega'$, so that $\Omega'\cap\Omega=\varnothing$.
Lemma~\ref{lem:block}, applied to $\Omega'$, yields signs $\varepsilon_k,\varepsilon'_k\in\{-1,1\}$,
a nonempty set $W\subset\Omega'$ and functions
$$
f_1=(3N)^{-\mu}N^{-1/s}\sum_{k\in\Omega'}\varepsilon_k h_k,
\qquad
f_2=(3N)^{-\mu}(\card W)^{-1/s}\sum_{k\in W}\varepsilon'_k h_k
$$
of the class $W_s^\mu$ whose coefficients vanish outside $\Omega'$. Since
$-\mu+3/4-1/s=-(\mu-a(s))$ and $-\mu+1/4-1/(3s)=-(\mu-b(s))$ by \eqref{exps}, the bounds of the same lemma
can be written as
$$
\|f_1\|_\infty\ge\tfrac1{10}\,3^{-\mu}N^{-(\mu-a(s))},
\qquad
\|f_2\|_\infty\ge c_*\,3^{-\mu}N^{-(\mu-b(s))} .
$$
Set $f=f_1$ for $s>4/3$ and $f=f_2$ for $s<4/3$. In the first case $\bar a(s)=a(s)$, in
the second $\bar a(s)=b(s)$; therefore in both cases the exponent equals
$-(\mu-\bar a(s))=-\beta$ (see~\eqref{beta}) and $\|f\|_\infty\ge c\,N^{-\beta}$. Since $G_{\Omega}f=\mathbf 0$, we have
$D_{\Omega}(\delta)\ge\|f\|_\infty$, and Lemma~\ref{lem:rd}(a) gives
\[e_\delta(W_s^\mu,\psi G_{\Omega},C(\mathbb R),\ell_p)\ge c\,N^{-\beta};\]
it remains to take
the infimum over $\Psi_N$. For $1\le N<N_0$, the inclusion $\Psi_N\subset\Psi_{N_0}$ implies
$$
R_{N,\delta}\ge R_{N_0,\delta}\ge c\,N_0^{-\beta}\ge c\,N_0^{-\beta}N^{-\beta}.
$$
Thus, $R_{N,\delta}\ge c\max\{\delta^e,N^{-\beta}\}\ge\frac c2\left(\delta^{e}+N^{-\beta}\right)$;
this is the left inequality of \eqref{eq:Rord}.

Let $\psi$ be a consistent recovery operator for $\Omega=\Omega_N$. Its values are continuous by
(I1), and $\psi$ uses only the coordinates with indices in $\Omega_N$, so that
$\psi G_{\Omega_N}\in\Psi_N$. On admissible data $\mathbf y=G_{\Omega_N}f+P_{\Omega_N}\bxi$
Lemma~\ref{lem:rd}(b) gives $\|f-\psi(\mathbf y)\|_\infty\le2D_{\Omega_N}(\delta)$. By
\eqref{reduction} and Theorem~\ref{thm:Eopt}, $D(\delta)\le\mathcal E_\delta\le C\delta^{e}$. Therefore,
for $N\ge2$, by Lemma~\ref{lem:rd}(c)
\[
\|f-\psi(\mathbf y)\|_\infty\le2D(\delta)+2cN^{-\beta}\le2C\delta^{e}+2cN^{-\beta},
\]
while for $N=1$, since $D_{\Omega_1}(\delta)\le B$ and $N^{-\beta}=1$, we have
\[
\|f-\psi(\mathbf y)\|_\infty\le2B=2BN^{-\beta}.
\]
Taking the supremum over $f\in W_s^\mu$ and $\|\bxi\|_{\ell_p}\le\delta$, we obtain for all $N\ge1$
\[
R_{N,\delta}(W_s^\mu,\Psi_N,C(\mathbb R),\ell_p)\le
e_\delta(W_s^\mu,\psi G_{\Omega_N},C(\mathbb R),\ell_p)\le C\left(\delta^{e}+N^{-\beta}\right),
\]
where the first inequality follows from $\psi G_{\Omega_N}\in\Psi_N$. This is the right inequality of
\eqref{eq:Rord}. $\Box$

In region III the two terms in \eqref{eq:Rord} reflect different obstructions.
The proof of the lower bound in $\delta$ uses the function $f_1$:
admissible noise hides its coefficients even under full information.
The lower bound in $N$ uses the function $f_2$, whose coefficients
lie outside the observed set; it remains indistinguishable from the opposite
function even under exact observations. The lower bound for $f_1$ contains a power of $N$ with
exponent $-(\mu-a(s))$, and that for $f_2$, with exponent $-(\mu-b(s))$. In region III
$s<4/3$, so that $b(s)>a(s)$ (Section~\ref{sec:exps}) and the bound for $f_2$ decreases in $N$
more slowly; it is this bound that gives the exponent $\beta=\mu-b(s)$. This explains the different roles of the two
constructions of Lemma~\ref{lem:block}.

\subsection{Explicit linear method in region III}
Recall that $M_0=\delta^{-1/(\mu+\lambda)}$,\, $M_1=\delta^{-1/(\mu+\lambda/3)}$,\,
$\lambda=a(p)-a(s)=1/s-1/p>0$, and $\Lambda(m)$ is defined before Lemma~\ref{lem:levels}.

Lemma~\ref{lem:levels} gives a suitable truncation level at each point,
but to construct a continuous function one must also pass continuously
between neighboring integer levels. Therefore, below, partial sums with two
neighboring upper limits are joined linearly. The factor $\eta(m)$
gradually reduces the result to zero between $M_1$ and $2M_1$:
on this interval the function itself is already small by~\eqref{eq:lin-P2}.
Finally, replacing $\Lambda$ by $\Lambda_N$ keeps only the first $N$
coefficients and leads to the additional term $N^{-\beta}$ in the error bound.

Let the parameters $p, s$ lie in region III and $m=x^2/2$. Set
$$
\eta(m)=\begin{cases}
1, & 0\le m\le M_1,\\
(2M_1-m)/M_1, & M_1<m<2M_1,\\
0, & m\ge2M_1.
\end{cases}
$$
For an integer $N\ge1$ define the capped level
$\Lambda_N(m)=\min\{\Lambda(m),N-1\}$ and at each point $x$ set
\[
L=\Lambda_N(m),\qquad n=\lfloor L\rfloor,\qquad \varphi=L-n,\qquad 0\le\varphi<1.
\]
Here $L$ is the truncation level at the point $x$: it is the level $\Lambda(m)$ chosen before
Lemma~\ref{lem:levels}, reduced if necessary to $N-1$, so that the method uses
at most $N$ coefficients. The level $L$ need not be an integer, while a sum can only run up to
an integer index. Therefore $n$ is the index of the last term that enters the sum in full, and
$\varphi$ is the fractional part of $L$, that is, the fraction with which the next, $(n+1)$st, term is taken.
If the sum stopped at the index $n$, then, when $L$ passes through an integer, the value of
the method would jump. The fraction $\varphi$ itself jumps when the level passes through an integer, but the sum
in which the next term is taken with this fraction varies continuously with $x$; this removes the jump.

The method takes the first $n+1$ data values with the factor $\eta(m)$, the next one with the fraction $\varphi$
of this factor, and does not use the rest:
\begin{equation}\label{eq:lin-method}
\psi_N(\mathbf y)(x)=\sum_{k=0}^{N-1}\omega_k^N(x)\,y_k h_k(x),
\qquad
\omega_k^N(x)=
\begin{cases}
\eta(m), & k\le n,\\
\eta(m)\,\varphi, & k=n+1,\\
0, & k\ge n+2.
\end{cases}
\end{equation}
In other words, the partial sum up to $n$ is supplemented by the fraction $\varphi$ of the next term and then
multiplied by $\eta(m)$:
\begin{equation}\label{psi_N}
\psi_N(\mathbf y)(x)=\eta(m)\bigl[(1-\varphi)S_n\mathbf y(x)
+\varphi S_{n+1}\mathbf y(x)\bigr],
\qquad
S_j\mathbf y=\sum_{k=0}^{j}y_k h_k .
\end{equation}
The same weights are written as a single formula in terms of the level:
\[
\omega_k^N(x)=\eta(m)\min\bigl\{1,\max\{0,\Lambda_N(m)+1-k\}\bigr\}.
\]
Indeed, $L+1-k\ge1+\varphi$ for $k\le n$, $L+1-k=\varphi$ for $k=n+1$, and
$L+1-k\le\varphi-1<0$ for $k\ge n+2$. This form is needed for the continuity of the method in $x$.

The number $N$ is the amount of information, common to all points $x$: the method uses only the coordinates
$y_0,\ldots,y_{N-1}$. At a given point $x$ only the indices $k\le n+1$ are actually used,
and the index $n+1$ only when $\varphi>0$. All of them are less than $N$. Indeed,
$n\le L\le N-1$. If $\varphi>0$, then $L$ is not an integer, while $N-1$ is an integer; hence $L<N-1$ and
$n+1\le N-1$. The number of terms actually summed at the point $x$ varies with $x$ together with
the level $L$, while the upper limit $N-1$ in \eqref{eq:lin-method} does not depend on $x$.

\begin{figure}[htb]
\centering
\includegraphics[width=\textwidth]{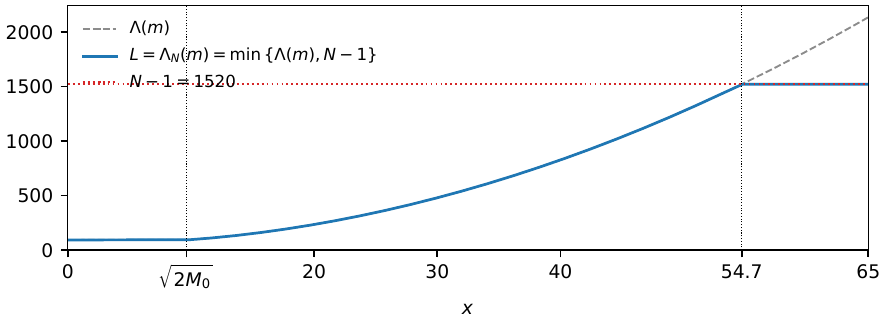}
\caption{The level $L=\Lambda_N(m)$, $m=x^2/2$, the uncapped level $\Lambda(m)$, and
the number $N-1$ for $s=1$, $p=2$, $\mu=0.1$, $\delta=0.1$, $N=1521$; here $M_0\approx46.4$, and
$L=2M_0$ for $x\le\sqrt{2M_0}\approx9.6$.}\label{fig:levels}
\end{figure}

Figure 2 shows that the number $L$ of terms actually summed coincides with $N-1$ only for $x\ge 54.7$.
The factor $\eta(m)$
differs from one only for $m>M_1\approx5623$, that is,
for $x>\sqrt{2M_1}\approx106$.

\begin{proposition}\label{prop:linear3}
For every integer $N\ge1$ the
following statements hold.
\begin{itemize}
\item[\rm(i)] The operator $\psi_N$ is linear, maps every real sequence to a
function continuous on $\mathbb R$ that vanishes for $x^2\ge4M_1$, and uses only the coordinates with
indices $k<N$. In particular, $\psi_N G_{\Omega_N}\in\Psi_N$.
\item[\rm(ii)] The following bound holds:
\begin{equation}\label{eq:lin-rate}
e_\delta(W_s^\mu,\psi_N G_{\Omega_N},C(\mathbb R),\ell_p)\le C\left(\delta^{e}+N^{-\beta}\right);
\end{equation}
therefore, by Theorem~\ref{thm:Rord}, the method $\psi_N G_{\Omega_N}$ is order-optimal in $\Psi_N$,
and for $N\ge\delta^{-e/\beta}$ its error is at most $C\delta^{e}$.
\item[\rm(iii)] The method $\psi_N$ is not diagonal. The error of every diagonal method on a
finite (possibly empty) index set is at least
$\min\{c\,\delta^{e_{\mathrm{tr}}},\pi^{-1/4}\}$ with the constant $c$ from (I3), and
$e_{\mathrm{tr}}<e$.
\end{itemize}
\end{proposition}

\textit{Proof.}
For $m=x^2/2\ge 2M_1$ we have $\eta(m)=0$, and hence, by \eqref{eq:lin-method}, the operator $\psi_N$ also equals $0$.
The continuity of $\Lambda$ follows from the equalities
$t_*(M_0)=M_0$, $t_*(M_1)=M_1^{1/3}$; the function $\eta$ is continuous as well.
Therefore, by the form of the weights in terms of the level, all $\omega_k^N$ are continuous, and \eqref{eq:lin-method} is a finite sum,
linear in the data and using only the first $k<N$ coordinates. For $\eta=0$ it vanishes.
This proves (i), including $\psi_N G_{\Omega_N}\in\Psi_N$.

\emph{Error bound.} Let $\mathbf y=\mathbf f+\bxi$ be admissible data. The
Fourier--Hermite series of $f$ converges at every point (I1); hence, for an integer $j\ge0$,
\[
f(x)-S_j\mathbf y(x)=\sum_{k>j}a_k(f)h_k(x)+\sum_{0\le k\le j}(-\xi_k)h_k(x).
\]
Here $\|(\underline k^{\,\mu}a_k(f))\|_{\ell_s}=\|f\|_{s,\mu}\le1$ and
$\|-\bxi\|_{\ell_p}\le\delta$. Therefore \eqref{split} with $L=j$ gives
\begin{equation}\label{eq:lin-split}
|f(x)-S_j\mathbf y(x)|\le\mathcal T(x,j)+\delta\,\mathcal H(x,j),
\qquad j\ge0 .
\end{equation}

Fix $x$ and use $L$, $n$, $\varphi$ from the definition of the method.
Denote by $Q$ the interpolated value from the representation of the method
by partial sums (see~\eqref{psi_N}):
\begin{equation}   \label{Q}
Q=(1-\varphi)S_n\mathbf y(x)+\varphi S_{n+1}\mathbf y(x),
\qquad \psi_N(\mathbf y)(x)=\eta(m)Q.
\end{equation}
We first estimate the error of this interpolation.

The sums in $\mathcal T$ and $\mathcal H$ are taken over integer $k$, and $n=\lfloor L\rfloor$; hence
$\mathcal T(x,n)=\mathcal T(x,L)$ and $\mathcal H(x,n)=\mathcal H(x,L)$; moreover,
$\mathcal T(x,n+1)\le\mathcal T(x,n)$, since $\mathcal T$ is nonincreasing in its second argument.
The triangle inequality and (E0) give $\mathcal H(x,n+1)\le\mathcal H(x,L)+C_0$.
Consequently, for $j=n,n+1$
\begin{equation}\label{eq:lin-E1}
\mathcal T(x,j)+\delta\,\mathcal H(x,j)\le\mathcal T(x,L)+\delta\,\mathcal H(x,L)+C_0\delta .
\end{equation}
Since $0\le\varphi<1$, inequalities \eqref{eq:lin-split} and \eqref{eq:lin-E1} give
\[
|f(x)-Q|\le(1-\varphi)|f(x)-S_n\mathbf y(x)|
+\varphi|f(x)-S_{n+1}\mathbf y(x)|
\le\mathcal T(x,L)+\delta\mathcal H(x,L)+C_0\delta.
\]
Now, taking (\ref{Q}) into account, we have
\[
f(x)-\psi_N(\mathbf y)(x)=\eta(m)\bigl(f(x)-Q\bigr)+(1-\eta(m))f(x).
\]
Here $0\le\eta(m)\le1$, and $|f(x)|\le\mathcal T(x,-1)$ by \eqref{split} with $L=-1$.
Consequently,
\begin{equation}   \label{ineq}
|f(x)-\psi_N(\mathbf y)(x)|
\le\eta(m)\bigl[\mathcal T(x,L)+\delta\mathcal H(x,L)+C_0\delta\bigr]
+(1-\eta(m))\mathcal T(x,-1).
\end{equation}

Let $\delta<\delta_0$. We estimate both terms on the right-hand side of \eqref{ineq}.
The first one uses the level $L=\min\{\Lambda(m),N-1\}$. It equals one of the numbers
$\Lambda(m)$, $N-1$ and is at most each of them, while $\mathcal T$ is nonincreasing and $\mathcal H$ is
nondecreasing in the second argument. Therefore
\[
\mathcal T(x,L)\le\mathcal T(x,\Lambda(m))+\mathcal T(x,N-1),
\qquad
\mathcal H(x,L)\le\mathcal H(x,\Lambda(m)).
\]
With a constant $C=C(s,\mu)$ we have
\begin{equation}\label{add}
\mathcal T(x,N-1)\le CN^{-\beta},\qquad N\ge1.
\end{equation}
This is the bound (I4).

The first term in \eqref{ineq} is nonzero only for $m<2M_1$. For such $m$ relation
\eqref{eq:lin-P1} bounds the series remainder and the noise,
$\mathcal T(x,\Lambda(m))+\delta\mathcal H(x,\Lambda(m))$, by $C\delta^e$. Moreover,
$C_0\delta\le C_0\delta^e$, since $e\le1$ (Section~\ref{sec:exps}) and $\delta<1$. Together with
\eqref{add} and $\eta(m)\le1$ this gives: the first term is at most $C(\delta^e+N^{-\beta})$.
The second term is nonzero only for $m>M_1$,
where \eqref{eq:lin-P2} gives $\mathcal T(x,-1)\le C\delta^e$.
This yields \eqref{eq:lin-rate} for $\delta<\delta_0$.

Now let $\delta_0\le\delta<1$. Inequality \eqref{ineq} holds for all $\delta$;
we show that its right-hand side is bounded by a constant. The quantity $\mathcal T$ is nonincreasing
in the second argument, and the term with $k=0$ in $\mathcal T(x,-1)$ equals $|h_0(x)|$, since
$\underline0=1$. Therefore the triangle inequality in $\ell_{s'}$, (E0), and \eqref{add} with $N=1$
give
\[
\mathcal T(x,L)\le\mathcal T(x,-1)\le|h_0(x)|+\mathcal T(x,0)\le C_0+C .
\]
The norm in $\ell_{p'}$ is at most the norm in $\ell_1$; hence, by (E0),
$\mathcal H(x,L)\le C_0(n+1)\le C_0(L+1)$. The first term in \eqref{ineq} is nonzero only
for $m<2M_1$, and there $L\le\Lambda(m)<4M_1$. This is seen from the three ranges in the definition of
$\Lambda$: $\Lambda(m)=2M_0<2M_1$ for $m\le M_0$;
$\Lambda(m)=m+t_*(m)\le2m\le2M_1$ for $M_0<m\le M_1$, since $t_*(m)\le m$ there;
$\Lambda(m)=m+m^{1/3}\le2m<4M_1$ for $M_1<m<2M_1$, since $m>M_1>1$.
Finally, $M_1\le\delta_0^{-1/(\mu+\lambda/3)}$. Therefore the right-hand side of \eqref{ineq} is at most
a constant $C$ independent of $N$, $x$, $f$ and $\bxi$. By the paragraph
``Large $\delta$'' in Section~\ref{sec:unlim}, it is at most $C\delta_0^{-e}\delta^e$,
and \eqref{eq:lin-rate} is proved on the whole interval $0<\delta<1$.
Comparison with the lower bound \eqref{eq:Rord} gives the order-optimality of the method
in $\Psi_N$. For $N\ge\delta^{-e/\beta}$ we have
$N^{-\beta}\le\delta^e$; hence the error is at most $C\delta^e$.
This proves (ii).

For (iii), take specific data: $y_0=1$ and $y_k=0$ for $k\ge1$.
Since $L=\Lambda_N(m)\ge0$, we have $n\ge0$, and by \eqref{eq:lin-method} the weight $\omega_0^N(x)$
equals $\eta(m)$. Therefore our method recovers the function $\eta(x^2/2)h_0(x)$.
At the origin it equals $h_0(0)\ne0$, and for $x^2\ge4M_1$ it vanishes.
The diagonal method $\psi_\nu$ from (I3) gives on the same data $\nu_0h_0$ if
$0\in\Omega$, and the zero function if $0\notin\Omega$. In both cases this is a function $c\,h_0$
with a factor $c$ independent of $x$ (for $0\notin\Omega$, $c=0$). If
$c\,h_0(x)=\eta(x^2/2)h_0(x)$ for all $x$, then $x=0$ would give $c=1$, since
$h_0(0)\ne0$, while $x^2\ge4M_1$ would give $c=0$, since $h_0$ vanishes nowhere.
Therefore no function of the form $c\,h_0$ coincides with the function obtained.
Consequently, $\psi_N$ is not diagonal. For every diagonal method
with a nonempty index set, the lower bound $c\delta^{e_{\mathrm{tr}}}$
is given in (I3). If the set is empty, the method uses no coefficient and, for
any data, gives the zero function. Take $f=h_0$. By the orthonormality of the Hermite
functions, its coefficients are $a_0(h_0)=1$ and $a_k(h_0)=0$ for $k\ge1$. Since
$\underline0=1$, we have $\|h_0\|_{s,\mu}=1$, that is, $h_0\in W_s^\mu$. Feed the method exact
data ($\bxi=\mathbf 0$ is an admissible error). The method returns the zero function, and its
error equals $\|h_0\|_\infty$. Since $h_0(x)=\pi^{-1/4}e^{-x^2/2}$, the maximum
is attained at the point $x=0$ and $\|h_0\|_\infty=\pi^{-1/4}$.
This gives the stated general lower bound. The inequality
$e_{\mathrm{tr}}<e$ follows from \eqref{eq:gap3}. $\Box$

\vskip 2mm

\begin{remark}\label{rem:arithmetic-cost}\rm
\textit{On the number of arithmetic operations.}
Let, as before, $L=\Lambda_N(m)$, $n=\lfloor L\rfloor$ and $\varphi=L-n$. By
\eqref{eq:lin-method}, the weights at the point $x$ are determined by three numbers $\eta(m)$, $n$, $\varphi$.
With $M_0$, $M_1$ computed in advance, these numbers are found in a bounded number of operations independent of
$N$, and, with the values $h_k(x)$ known, computing $\psi_N(\mathbf y)(x)$
requires $O(n+2)$ arithmetic operations.
The simplicity of the explicit weights is combined with the linearity of the method, the continuity of the recovered
function, and order-optimality for every fixed number $N$ of used coefficients in region III
(Proposition~\ref{prop:linear3}).

\end{remark}

\begin{remark}\label{rem:summation-comparison}\rm
\textit{On the relation to other summation methods.}
The choice of a point-dependent order in spectral recovery has appeared before.
For example, in \cite[pp.~156--157, Theorem 3.1]{TadmorTanner2002}, to recover
piecewise smooth functions from spectral data and to suppress the Gibbs phenomenon
away from discontinuities, a localizing kernel is used that contains the Dirichlet kernel
of degree $p(x)\asymp d(x)N$, where $d(x)$ is the distance to the nearest discontinuity.
In \cite[\S~1, Theorem 2.1]{TadmorTanner2005}
the order of smoothing is also chosen depending on the point,
while the common limit $N$ in the Fourier sum is kept. These works give examples
of choosing point-dependent parameters in recovery; here, a point-dependent
truncation level is used for minimax recovery
from inexact Fourier--Hermite coefficients.
\end{remark}

\vskip 2mm

\subsection{Information complexity}

It remains to find the smallest number of coefficients for which the restriction
on the amount of information no longer changes the order of the error. By
Theorem~\ref{thm:Rord}, for this it is necessary and sufficient to make
$N^{-\beta}$ and $\delta^e$ of the same order. We state the result taking into account the constant $K$
from definition~\eqref{eq:Nmin}.

\begin{theorem}\label{thm:complexity}
There exists $K_0=K_0(s,p,\mu)$ such that for every $K\ge K_0$ the set on the right-hand side of
\eqref{eq:Nmin} is nonempty and
\begin{equation}\label{eq:Nminord}
N_{\min}(\delta)\asymp\delta^{-e/\beta}
\end{equation}
with constants depending on $s$, $p$, $\mu$ and $K$. The difference
$e/\beta-1/\rho=(e-e_{\mathrm{tr}})/\beta$, by \eqref{eq:gap3}, equals zero outside region III, where
$N_{\min}(\delta)\asymp\delta^{-1/\rho}$ is the order of balanced truncation (I2), and is positive in
region III\@. For every $K>0$, the inequality $R_{N,\delta}\le K\delta^{e}$ implies
\[N\ge(c/K)^{1/\beta}\delta^{-e/\beta}\ge(c/K)^{1/\beta}\delta^{-1/\rho};\]
below the balanced level, the order $\delta^{e}$ of the radius is not attained.
\end{theorem}

\textit{Proof} follows the scheme of \cite[Theorem~4.12]{SolodkyVolynets2026}. For
$N_1=\lceil\delta^{-e/\beta}\rceil$ we have $N_1^{-\beta}\le\delta^{e}$, and by \eqref{eq:Rord}
$$
R_{N_1,\delta}(W_s^\mu,\Psi_{N_1},C(\mathbb R),\ell_p)\le2C\,\delta^{e} ;
$$
hence, for $K\ge K_0=2C$ the set in \eqref{eq:Nmin} is nonempty and
$$
N_{\min}(\delta)\le N_1\le2\delta^{-e/\beta}.
$$
Conversely, for every $K>0$, if $R_{N,\delta}\le K\delta^{e}$, then by \eqref{eq:Rord}
$$
c\,N^{-\beta}\le K\delta^{e},
\qquad\text{whence}\qquad
N\ge\Bigl(\frac cK\Bigr)^{1/\beta}\delta^{-e/\beta},
$$
and since $e/\beta\ge1/\rho$ and $\delta<1$, this is at least $(c/K)^{1/\beta}\delta^{-1/\rho}$; for
$K\ge K_0$ this also holds for $N=N_{\min}(\delta)$. $\Box$

In region III, for $K\ge K_0$ and $N\ge N_{\min}(\delta)$ we have $R_{N,\delta}\le K\delta^{e}$, since
$R_{N,\delta}$ is nonincreasing in $N$; therefore, by Proposition~\ref{prop:linear3}\,(iii),
the ratio of the error of every diagonal method on at most $N$ coefficients to
$R_{N,\delta}$ is at least $c_K\delta^{-(e-e_{\mathrm{tr}})}\to\infty$, where
$c_K=c_K(s,p,\mu,K)>0$, and no order-optimal method is diagonal.

\begin{remark}\label{rem:ladder}\rm
For $s<4/3<p$ and $\mu>\max\{0,11/12-1/s\}$ we have $\tilde e<e_{\mathrm{tr}}<e$, where
$\tilde e=(\mu+1/s-11/12)/(\mu+1/s-1/p)$ is the exponent of the standard scheme that bounds
each $h_k$ by its maximum before summation (\cite[Theorem~3.1]{MathePereverzev2002},
\cite[Remark~4.8]{HerMinCompanion}); for example, for $s=1$, $p=2$, $\mu=1$ this is
$13/18<13/16<5/6$. The gain of $e$ over $e_{\mathrm{tr}}$ comes from choosing the truncation level separately at
each point $x$.
\end{remark}

\begin{remark}\label{rem:wellposed}\rm
Following \cite[\S\,2.8]{LuP}, we call a problem well-posed if its solution can be approximated
with accuracy of order $O(\delta)$. By \eqref{beta}, the exponent $e$ equals one exactly
when $a_+(p)=0$, that is, for $p<4/3$. Therefore Theorem~\ref{thm:Eopt} gives the following.
For $p<4/3$ (regions II and IV) the problem of summation of the Fourier--Hermite series from perturbed
coefficients in the metric $C(\mathbb R)$ is well-posed:
\[
c\,\delta\le\mathcal E_\delta(W_s^\mu,\Psi,C(\mathbb R),\ell_p)\le C\,\delta .
\]
For $p>4/3$ (regions I and III) we have $e<1$, so that
$\mathcal E_\delta(W_s^\mu,\Psi,C(\mathbb R),\ell_p)/\delta\ge c\,\delta^{e-1}\to\infty$
as $\delta\to0$: the order $\delta$ is attained by no recovery operator, and the problem
is ill-posed. In region III this was previously established only for finite diagonal
methods \cite[Remark~5.2]{HerMinCompanion}.
\end{remark}

\section{Numerical example}\label{sec:numerics}

\emph{Recovering a wave function from the amplitudes of the oscillator levels.} In dimensionless units ($\hbar=m=\omega=1$) the Hamiltonian of the harmonic
oscillator is $\frac12\bigl(-\frac{d^2}{dx^2}+x^2\bigr)$, and the Hermite functions $h_k$ are
its eigenfunctions with energy levels $k+\frac12$. A particle in the attractive
potential
\[
V(x)=-\frac1{1+|x|}
\]
has a bound state with energy $E=-\frac12$: the equation $-\frac12f''+Vf=Ef$ is solved by
the function normalized in ${\cal L}_2(\mathbb R)$
\[
f_0(x)=\sqrt{2/5}\,(1+|x|)\,e^{-|x|}.
\]
Let the potential $V$ be replaced by the oscillator potential $x^2/2$ within
a time that is short compared with the periods of transitions from the state $f_0$ to the oscillator levels.
Then during the replacement
the wave function does not change, and the particle ends up on level $k$ of the oscillator with
probability $a_k(f_0)^2$ \cite[\S\,41, (41.6)]{LandauLifshitz}. The amplitudes $a_k(f_0)$, together with
their signs, found with an error, are exactly the data of the summation problem: from them one has to
recover the wave function on the whole axis in the uniform metric. The data error is measured in $\ell_2$,
that is, $p=2$.

\emph{Coefficients and class.} The function $f_0$ is even, so $a_k(f_0)=0$ for odd $k$, and
for even indices $a_k(f_0)=2\sqrt{2/5}\,(I_k+J_k)$, where
\[
I_n=\int_0^\infty e^{-x}h_n(x)\,dx,\qquad J_n=\int_0^\infty x\,e^{-x}h_n(x)\,dx .
\]
Integration by parts and the formula $h_n'=\sqrt{n/2}\,h_{n-1}-\sqrt{(n+1)/2}\,h_{n+1}$ give
\[
\sqrt{\tfrac{n+1}{2}}\,I_{n+1}=\sqrt{\tfrac n2}\,I_{n-1}+h_n(0)-I_n,
\qquad
\sqrt{\tfrac{n+1}{2}}\,J_{n+1}=\sqrt{\tfrac n2}\,J_{n-1}+I_n-J_n,
\qquad n\ge0,
\]
where $I_{-1}=J_{-1}=0$. By these formulas the coefficients for $k\le338764$ were computed in
arithmetic with 387 decimal digits; the sum $\sum_{k\le6\cdot10^4}a_k(f_0)^2$ agrees with
$\|f_0\|^2_{{\cal L}_2(\mathbb R)}=1$ to within $10^{-15}$. Near zero
$f_0(x)=\sqrt{2/5}\,\bigl(1-\frac{x^2}2+\frac{|x|^3}3+\ldots\bigr)$, therefore
\[
|a_k(f_0)|\approx0.424\,k^{-9/4},\qquad k\ \text{even},
\]
and $\|f_0\|_{1,\mu}<\infty$ for $\mu<5/4$. Take $s=1$, $\mu=1$ and
$f=f_0/R$, where $R=\|f_0\|_{1,\mu}\approx5.087$, so that $\|f\|_{1,\mu}=1$; the data are
$y_k=a_k(f)+\xi_k$, $\|\bxi\|_{\ell_2}\le\delta$. The pair $s=1$, $p=2$
lies in region III, and by \eqref{beta} $\beta\approx1.083$, $\rho\approx1.333$,
$e_{\mathrm{tr}}=0.8125$, $e\approx0.833$.

\emph{Methods.} Balanced truncation $S_{N-1}\mathbf y$ from (I2) uses
$N=\lceil\delta^{-1/\rho}\rceil+1$ coefficients, and the method $\psi_N$ from
\eqref{eq:lin-method} uses $N=\lceil\delta^{-e/\beta}\rceil+1$ coefficients
(Proposition~\ref{prop:linear3}(ii)).

\emph{Noise.} The noise is random: $\bxi=\delta\mathbf v/\|\mathbf v\|_{\ell_2}$, where
$v_0,\ldots,v_{N-1}$ are independent and standard normal. For each $\delta$ and each
method $100$ realizations are taken; the noise samples of the two methods are different, since their values of $N$ differ.
The supremum over $x$ is taken on the interval $|x|\le\sqrt{2N+1}+6$ on a grid with step
$\min\{0.0025,\,2\pi/(10\sqrt{2N+1})\}$, that is, with at least ten points per period
$2\pi/\sqrt{2N+1}$;
the values $h_k(x)$ are computed by the three-term recurrence relation.

\begin{table}[ht]
\centering
\caption{Error $\sup_x|f(x)-\psi(\mathbf y)(x)|$, mean over $100$ realizations
of random noise of norm $\delta$; $N$ is the number of coefficients used.}\label{tab:random}
\begin{tabular}{llllll}
\hline
$\delta$ & $N$ (truncation) & $N$ (method $\psi_N$) & truncation $S_{N-1}$ & method $\psi_N$ & ratio\\
\hline
0.1    & 7     & 7     & $7.91\cdot10^{-2}$ & $7.81\cdot10^{-2}$ & 1.01\\
0.05   & 11    & 12    & $3.79\cdot10^{-2}$ & $3.79\cdot10^{-2}$ & 1.00\\
0.01   & 33    & 36    & $6.86\cdot10^{-3}$ & $6.77\cdot10^{-3}$ & 1.01\\
0.005  & 55    & 60    & $3.27\cdot10^{-3}$ & $3.20\cdot10^{-3}$ & 1.02\\
0.001  & 179   & 205   & $5.67\cdot10^{-4}$ & $5.34\cdot10^{-4}$ & 1.06\\
0.0005 & 301   & 348   & $2.54\cdot10^{-4}$ & $2.41\cdot10^{-4}$ & 1.06\\
0.0001 & 1001  & 1195  & $4.21\cdot10^{-5}$ & $3.73\cdot10^{-5}$ & 1.13\\
0.00005 & 1683 & 2036  & $1.91\cdot10^{-5}$ & $1.68\cdot10^{-5}$ & 1.14\\
0.00001 & 5625 & 7019  & $3.10\cdot10^{-6}$ & $2.56\cdot10^{-6}$ & 1.21\\
0.000005 & 9459 & 11961 & $1.39\cdot10^{-6}$ & $1.13\cdot10^{-6}$ & 1.23\\
0.000001 & 31624 & 41248 & $2.21\cdot10^{-7}$ & $1.67\cdot10^{-7}$ & 1.32\\
0.0000005 & 53184 & 70300 & $9.87\cdot10^{-8}$ & $7.36\cdot10^{-8}$ & 1.34\\
0.0000001 & 177829 & 242448 & $1.55\cdot10^{-8}$ & $1.09\cdot10^{-8}$ & 1.42\\
\hline
\end{tabular}
\end{table}

\begin{figure}[ht]
\centering
\includegraphics[width=0.8\textwidth]{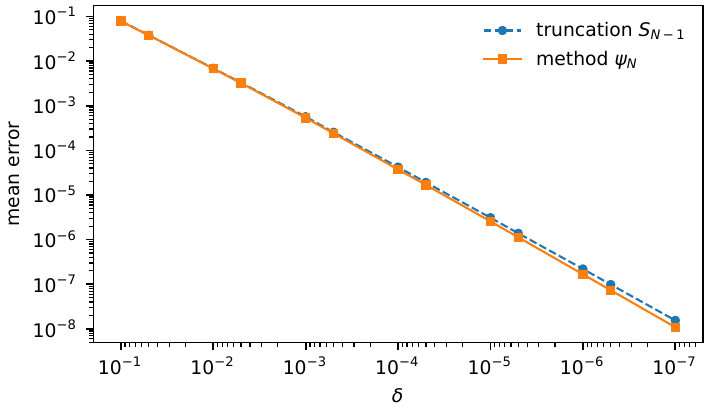}
\caption{Mean error of the truncation $S_{N-1}$ and of the method $\psi_N$ as a function of $\delta$
(Table~\ref{tab:random}).}\label{fig:numerics}
\end{figure}

\begin{figure}[ht]
\centering
\includegraphics[width=\textwidth]{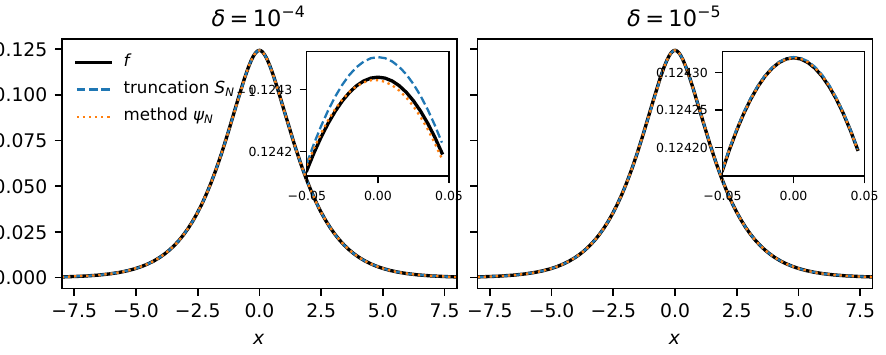}
\caption{The function $f$ and its approximations by the truncation $S_{N-1}$ and by the method $\psi_N$ from data with
a single realization of random noise for $\delta=10^{-4}$ and $\delta=10^{-5}$; the insets show
the interval $|x|\le0.05$.}\label{fig:approx}
\end{figure}

\emph{Results.} For $\delta\ge0.005$ both methods use almost the same number of
coefficients; the mean errors differ by a factor of at most $1.02$, and this difference
does not exceed the statistical error of the mean over $100$ realizations
(Table~\ref{tab:random}, Figure~\ref{fig:numerics}). For $\delta\le10^{-3}$ the method $\psi_N$
is on average more accurate than the truncation, and the gain grows: the ratio equals $1.06$ for $\delta=10^{-3}$ and
$1.42$ for $\delta=10^{-7}$. The gain is related to the noise: for $\delta\le10^{-4}$, near $x=0$
the method $\psi_N$ sums the coefficients only up to the index $2M_0$, which is smaller than $N-1$ for the truncation, and
hence sums fewer noisy terms. For $\delta=10^{-4}$ and $\delta=10^{-5}$
the graphs of both approximations at the scale of Figure~\ref{fig:approx} coincide with the graph of $f$.

We add that Table~\ref{tab:random} gives the mean error for a single
function $f$ under random noise.
At the same time, Theorem~\ref{thm:Eopt} and Proposition~\ref{prop:linear3}
estimate the error $e_\delta(W_1^\mu,\psi G,C(\mathbb R),\ell_2)$ on the whole class.
Our additional computations show that on the class $W_1^\mu$ for $p=2$
the guaranteed error of the truncation method is larger than that of the method $\psi_N$,
and the ratio of these quantities grows as $\delta$ decreases: from $1.04$ to $1.16$ for $\mu=0.1$
and $\delta=0.1,\,0.05,\,0.02$, and from $1.00$ to $1.12$ for $\mu=0.2$ and $\delta=0.1,\,0.05,\,0.01,\,0.005$.

\appendix

\section{Proof of Lemma~\ref{lem:levels}}\label{app:levels}
\renewcommand{\theequation}{A.\arabic{equation}}

\textit{Proof.} We prove the estimates in three steps: first we estimate the sums of powers of $|h_k(x)|$,
then the quantities $\mathcal H$, $\mathcal T$, and finally we substitute the chosen level.
From $\delta<\delta_0$ it follows that $M_0>8$. Next, $M_0<M_1$: by \eqref{M0} these are powers of the number
$\delta^{-1}>1$ with exponents $1/(\mu+\lambda)<1/(\mu+\lambda/3)$. First we assume $s>1$;
the case $s=1$ is treated at the end of Step~2.

\emph{Step 1. Estimates of sums over indices.} Let $m\ge8$, $0\le t\le m$, $L=m+t$, $\tau=m^{1/3}+t$ and $1\le r<\infty$. Then
\begin{equation}\label{K0}
\begin{gathered}
\sum_{k<m}|h_k(x)|^r\le C(m^{1/3})^{1-r/4},\\
\sum_{m\le k\le L}|h_k(x)|^r\le C\tau^{1-r/4}\quad(r<4),\\
\sum_{k>L}|h_k(x)|^r\le C\tau^{1-r/4}\quad(r>4),
\end{gathered}
\end{equation}
where $C$ depends only on $r$. The number $\tau=m^{1/3}+(L-m)$ measures the level $L$ from $m$, and
not from zero: by (E1), for $k\ge m$ the quantity $|h_k(x)|$ is estimated through $m^{1/3}+k-m$,
therefore the sums over $m\le k\le L$ and over $k>L$ are estimated through $\tau$, and not through $L$.
For the first estimate, we split the indices $k<m$ into the groups
$j\le(m-k)/m^{1/3}<j+1$, $j=0,1,\ldots$. Each group contains at most $m^{1/3}+1\le2m^{1/3}$
indices; by (E3) its contribution is at most $Cm^{-r/12}m^{1/3}e^{-crj^{3/2}}$.
The sum over $j$ is bounded by a constant: for $j\ge1$ we have $j^{3/2}\ge j$ and obtain
a convergent geometric progression. This gives the first estimate in \eqref{K0}.

For the other two estimates we set $v_k=m^{1/3}+k-m$; by (E1) we have
$|h_k(x)|^r\le Cv_k^{-r/4}$. The numbers $v_k$ are spaced by~1, and for $m\le k\le L$ they lie
between $m^{1/3}\ge2$ and $\tau$. We split this range into the groups
$\tau/2^{j+1}<v_k\le\tau/2^j$, $j\ge0$. A nonempty group contains at most
$\tau/2^j$ numbers, each of which contributes at most $C(\tau/2^j)^{-r/4}$.
Therefore, for $r<4$
\[
\sum_{m\le k\le L}|h_k(x)|^r
\le C\tau^{1-r/4}\sum_{j\ge0}2^{-j(1-r/4)}\le C\tau^{1-r/4}.
\]
For $k>L$ we have $v_k>\tau$. Now we take the groups $2^j\tau<v_k\le2^{j+1}\tau$:
each contains at most $2\cdot2^j\tau$ numbers, and each term is at most $C(2^j\tau)^{-r/4}$.
Consequently, for $r>4$
\[
\sum_{k>L}|h_k(x)|^r
\le C\tau^{1-r/4}\sum_{j\ge0}2^{-j(r/4-1)}\le C\tau^{1-r/4}.
\]
It is precisely the convergence of these two progressions that explains the different restrictions on $r$ in \eqref{K0}.

\emph{Step 2. Estimates of the noise and of the remainder of the series.} We estimate both quantities on the right-hand side of
\eqref{holder}: $\mathcal H$ collects the indices $0\le k\le L$ (the influence of the data error),
and $\mathcal T$ collects the indices $k>L$ (the remainder of the series).
In Step~3 the level is chosen differently in three ranges of $m$, and each of them needs its own estimates:
\begin{caselist}
\item[$m\le M_0$:] the level $L=2M_0\ge2m$ lies to the right of the transition zone; estimates
through $L$ itself, valid for all $x$, suffice (parts 2a, 2b);
\item[$M_0<m<2M_1$:] the level $L=m+t$, $0\le t\le m$, lies near the transition zone
$k\approx m$; one needs estimates through $\tau=m^{1/3}+t$, which can be much smaller than
$L$ (parts 2a, 2c);
\item[$m>M_1$:] the whole series is estimated, that is, $\mathcal T(x,-1)$ (part 2d).
\end{caselist}
Table~\ref{tab:step2} summarizes the levels and the estimates of Step~2 over these ranges.

\begin{table}[htb]
\centering
\caption{The level $L$, the number $\tau$, and the estimates of the remainder of the series $\mathcal T$ and of the quantity
$\mathcal H$ in the three ranges of $m$.}\label{tab:step2}
\renewcommand{\arraystretch}{1.3}
\small
\begin{tabular}{lllll}
\hline
range of $m$ & level $L$ & $\tau$ & $\mathcal T(x,L)$ & $\mathcal H(x,L)$\\
\hline
$m\le M_0$ & $L=2M_0\ge\max\{2m,1\}$ & --- & $\le CL^{a(s)-\mu}$, \eqref{K3} & $\le CL^{a(p)}$, \eqref{K1}\\
$M_0<m<2M_1$ & $L=m+t$, $0\le t\le m$ & $m^{1/3}+t$ & $\le Cm^{-\mu}\tau^{a(s)}$, \eqref{K4} & $\le C\tau^{a(p)}$, \eqref{K2}\\
$m>M_1$ & $L=-1$ (the whole series) & --- & $\le Cm^{b(s)-\mu}$, \eqref{K5} & $=0$\\
\hline
\end{tabular}
\end{table}
The case $s=1$ is treated in part 2e. Figure~\ref{fig:bounds} shows where
these estimates come from: for fixed $x$ the values $|h_k(x)|$ for $k<m$ are exponentially
small (E3), while for $k\ge m$ they decrease slowly under the curve $C(m^{1/3}+k-m)^{-1/4}$ (E1).
The sums of powers of $|h_k(x)|$ over the three zones $k<m$, $m\le k\le L$ and $k>L$ are estimated in~\eqref{K0}.

\emph{2a. Estimates of $\mathcal H$} are needed in Step~3 for $m\le M_0$ and for $M_0<m<2M_1$.
By \cite[Theorem~4.2]{HerMinCompanion} with $q=p'$ (there $1\le q<4$, here $1\le p'<4$),
for all $x$ and $N\in\mathbb N_0$
\[
\sum_{k=0}^{N}|h_k(x)|^{p'}\le C(N+1)^{1-p'/4}.
\]
We take $N=\lfloor L\rfloor$; since $\lfloor L\rfloor+1\le2L$ for $L\ge1$ and $1-p'/4>0$,
we get $\mathcal H(x,L)^{p'}\le CL^{1-p'/4}$. We take the $p'$-th root; since
\[
\frac{1-p'/4}{p'}=\frac1{p'}-\frac14=\frac34-\frac1p=a(p),
\]
we obtain
\begin{equation}\label{K1}
\mathcal H(x,L)\le CL^{a(p)},
\qquad x\in\mathbb R,\quad L\ge1 .
\end{equation}
Estimate \eqref{K1} does not take into account that for $k<m$ almost all $|h_k(x)|$ are exponentially small.
Under the conditions of Step~1 this fact allows us to refine \eqref{K1}. We split the sum $\mathcal H(x,L)^{p'}$
into the terms with $k<m$
and with $m\le k\le L$ and apply \eqref{K0} with $r=p'<4$: the first part is at most
$C(m^{1/3})^{1-p'/4}$, the second is at most $C\tau^{1-p'/4}$. Since $1-p'/4>0$ and
$m^{1/3}\le\tau$, both parts are at most $C\tau^{1-p'/4}$, and taking the $p'$-th root gives
\begin{equation}\label{K2}
\mathcal H(x,L)\le C\tau^{a(p)}.
\end{equation}

\emph{2b. The remainder of the series for a level far beyond the transition zone} is needed in Step~3 for
$m\le M_0$. Let $L\ge\max\{2m,1\}$. For $k>L$ we have $k\ge2m$ and $k\ge1$, so that
$\underline k=k$, and by (E2) $|h_k(x)|\le Ck^{-1/4}$. Therefore
\[
\mathcal T(x,L)^{s'}\le C\sum_{k>L}k^{-\sigma},
\qquad \sigma=s'\Bigl(\mu+\frac14\Bigr)>\frac{s'}4>1.
\]
The group $2^{j}L<k\le2^{j+1}L$, $j\ge0$, contains at most $2^{j}L+1$ indices; since $L\ge1$, this is
at most $2^{j+1}L$. Each term is at most $(2^{j}L)^{-\sigma}$. Since $\sigma>1$,
the progression over $j$ converges:
\[
\sum_{k>L}k^{-\sigma}
\le2L^{1-\sigma}\sum_{j\ge0}2^{-j(\sigma-1)}
\le CL^{1-\sigma}.
\]
Taking the $s'$-th root and using the equality $1/s'-1/4=a(s)$ give
\begin{equation}\label{K3}
\mathcal T(x,L)\le CL^{a(s)-\mu},
\qquad L\ge\max\{2m,1\} .
\end{equation}

\emph{2c. The remainder of the series for a level near the transition zone} is needed in Step~3 for
$M_0<m<2M_1$; there the level $L=m+t$ may be smaller than $2m$, and \eqref{K3} is not applicable.
Let the conditions of Step~1 hold. For $k>L\ge m$ we have $\underline k=k\ge m$, and since
$\mu>0$, we get $\underline k^{\,-\mu}\le m^{-\mu}$. We take this factor out and apply
\eqref{K0} with $r=s'>4$ to the zone $k>L$:
\[
\mathcal T(x,L)^{s'}\le m^{-\mu s'}\sum_{k>L}|h_k(x)|^{s'}\le Cm^{-\mu s'}\tau^{1-s'/4}.
\]
Since $(1-s'/4)/s'=a(s)$, taking the $s'$-th root gives
\begin{equation}\label{K4}
\mathcal T(x,L)\le Cm^{-\mu}\tau^{a(s)}.
\end{equation}

\emph{2d. The whole series} is estimated for Step~3 for $m>M_1$. Let $m\ge8$ and
$N=\lceil m/2\rceil$. By (I4), the indices $k\ge N$ contribute to $\mathcal T(x,-1)$ at most
$CN^{-\beta}\le Cm^{b(s)-\mu}$, since in region III $\bar a(s)=b(s)$ and $N\ge m/2$.
For $k<m/2$, by (E3), $|h_k(x)|\le Ce^{-cm}$ and $\underline k^{\,-\mu}\le1$; there are at most
$m$ such indices, so their part is at most $Cme^{-cm}\le Cm^{b(s)-\mu}$. The triangle
inequality in $\ell_{s'}$ gives
\begin{equation}\label{K5}
\mathcal T(x,-1)\le Cm^{b(s)-\mu} .
\end{equation}

\emph{2e. The case $s=1$.} Then $s'=\infty$ and
$\mathcal T(x,L)=\sup_{k>L}\underline k^{-\mu}|h_k(x)|$: instead of sums, the largest
term is estimated, directly by Lemma~\ref{lem:cap}. For \eqref{K3}, for $k>L\ge\max\{2m,1\}$ the estimate
(E2) gives $\underline k^{-\mu}|h_k(x)|\le Ck^{-\mu-1/4}\le CL^{-\mu-1/4}$. For \eqref{K4},
for $k>L=m+t$ we have $k-m>t$, and (E1) gives $\underline k^{-\mu}|h_k(x)|\le Cm^{-\mu}\tau^{-1/4}$.
Part 2d does not use the condition $s>1$: (I4) holds also for $s=1$, and the triangle inequality holds also in
$\ell_\infty$; therefore \eqref{K5} holds also for $s=1$. Since $a(1)=-1/4$ and
$b(1)=-1/12$, the estimates \eqref{K3}--\eqref{K5} hold also for $s=1$; the estimates \eqref{K1},
\eqref{K2} do not depend on $s$.

\emph{Step 3. Choice of the level in each range of $m$.} By \eqref{M0} the function $t_*(m)$
decreases in $m$, since its exponent $-\mu/\lambda$ is negative. We estimate the error at the
right end $m=M_1$ of the middle range. There $t_*(M_1)=M_1^{1/3}$, while $b(s)=a(s)/3$ and
$b(p)=a(p)/3$; therefore the equality $m^{-\mu}t_*^{a(s)}=\delta t_*^{a(p)}$, which defines
$t_*(m)$, takes at $m=M_1$ the form
\[
M_1^{b(s)-\mu}=M_1^{-\mu}t_*(M_1)^{a(s)}=\delta\,t_*(M_1)^{a(p)}=\delta M_1^{b(p)}.
\]
Since $t_*$ decreases, $M_0<M_1$ and $a(p)>0$, taking \eqref{M0e} into account we obtain
\[
\delta\,t_*(M_1)^{a(p)}\le\delta\,t_*(M_0)^{a(p)}=\delta M_0^{a(p)}=\delta^e .
\]
Thus,
\begin{equation}\label{eq:M1edge}
M_1^{b(s)-\mu}=\delta M_1^{b(p)}\le\delta^{e} .
\end{equation}
For $m\le M_0$ we use the general estimates. Since $M_0\ge m$ and $M_0>8$, we have
$\Lambda(m)=2M_0\ge\max\{2m,1\}$, and by \eqref{K1}, \eqref{K3} and \eqref{M0e}
\[
\mathcal T(x,\Lambda(m))+\delta\mathcal H(x,\Lambda(m))\le C(M_0^{a(s)-\mu}+\delta M_0^{a(p)})=2C\delta^e .
\]
For $M_0<m\le M_1$ we balance the estimates of the remainder of the series and of the noise. Here $m^{1/3}\le t_*(m)\le m$;
we apply \eqref{K2}, \eqref{K4} with $t=t_*(m)$, so that $t_*(m)\le\tau\le2t_*(m)$. Since
$a(s)<0<a(p)$, we have $\tau^{a(s)}\le t_*(m)^{a(s)}$ and $\tau^{a(p)}\le2^{a(p)}t_*(m)^{a(p)}$. The equality
defining $t_*(m)$ and the decrease of $t_*$ give
\[
\begin{aligned}
\mathcal T(x,\Lambda(m))+\delta\mathcal H(x,\Lambda(m))
&\le C\bigl(m^{-\mu}t_*(m)^{a(s)}+\delta t_*(m)^{a(p)}\bigr)\\
&=2C\delta\,t_*(m)^{a(p)}\le2C\delta\,t_*(M_0)^{a(p)}=2C\delta^e.
\end{aligned}
\]
For $M_1<m<2M_1$ the balance point already lies below $m^{1/3}$; we use the level
$\Lambda(m)=m+m^{1/3}$ and apply \eqref{K2}, \eqref{K4} with $t=m^{1/3}$. We use three
facts: $\tau=2m^{1/3}$, so that $\tau^{a(s)}\le Cm^{b(s)}$ and $\tau^{a(p)}\le Cm^{b(p)}$;
the exponent $b(s)-\mu$ is negative, and $m>M_1$; the exponent $b(p)$ is positive, and $m<2M_1$.
Together with \eqref{eq:M1edge} we obtain
\[
\mathcal T(x,\Lambda(m))+\delta\mathcal H(x,\Lambda(m))
\le C\bigl(m^{b(s)-\mu}+\delta m^{b(p)}\bigr)
\le C\bigl(M_1^{b(s)-\mu}+\delta M_1^{b(p)}\bigr)\le C\delta^e.
\]
Thus estimate \eqref{eq:lin-P1} is obtained on the whole range $m<2M_1$.
Finally, for $m>M_1$ the whole series is small by itself: by \eqref{K5}, since $b(s)-\mu<0$,
\[
\mathcal T(x,-1)\le Cm^{b(s)-\mu}\le CM_1^{b(s)-\mu}\le C\delta^e
\]
by \eqref{eq:M1edge}. This is \eqref{eq:lin-P2}. $\Box$

\end{document}